\documentclass[10pt, a4paper]{article}
\usepackage[margin=1in]{geometry}
\usepackage{float}

\usepackage{graphicx}
\usepackage{amsthm, amsfonts,amssymb} %amsmath,
\usepackage{lmodern}
\usepackage[T1]{fontenc}
\usepackage{mathtools}

\newtheorem{thm}{Theorem}[section]

\newtheorem{cor}[thm]{Corollary}

\newtheorem{mydef}{Definition}[section]

\newtheorem{rem}{Remark}[section]
\theoremstyle{remark}
\newcommand\bb[1]{\mathbf{#1}}
\newcommand\ddfrac[2]{\frac{\displaystyle #1}{\displaystyle #2}}

\newcommand\bint[1]{\displaystyle\int #1}

\usepackage[font=footnotesize,labelfont=bf]{caption}

\begin{document}
\date{}
\author{ Sr\dj{}an Trifunovi\'c \thanks{Department of Mathematics and Informatics, University of Novi Sad, Trg D. Obradovi\'{c}a 4, 21000 Novi Sad, Serbia,
    email: sergej1922@gmail.com}}

\title{Compressible fluids interacting with plates -- regularity and weak-strong uniqueness }

\maketitle

\begin{abstract}
	In this paper, we study a nonlinear interaction problem between compressible viscous fluids and plates. For this problem, we introduce relative entropy and relative energy inequality for the finite energy weak solutions (FEWS). First, we prove that for all FEWS, the structure displacement enjoys improved regularity by utilizing the dissipation effects of the fluid onto the structure and that all FEWS satisfy the relative energy inequality. Then, we show that all FEWS enjoy the weak-strong uniqueness property, thus extending the classical result for compressible Navier-Stokes system to this fluid-structure interaction problem.
\end{abstract}
 ${}$
\\${}$ 
% 76N10 :Existence, uniqueness, and regularity theory for compressible fluids and gas dynamics
% 35A02: Uniqueness problems for PDEs: global uniqueness, local uniqueness, non-uniqueness,
% 35R37: Moving boundary problems for PDEs
% 74F10: Fluid-solid interactions (including aero- and hydro-elasticity, porosity, etc.)
\section{Introduction}
The well-posedness theory for the nonlinear interaction problems between fluids and thin elastic structures (plates or shells) has seen quite a development in the last two decades. Such problems arise from various physical phenomena and they are quite challenging for various reasons. For example, the fluid domain changes in time and it depends on the displacement of the elastic structure - an unknown in the system, it is a mixed-type coupled physical system, usually there are both hyperbolic and parabolic equations in the system etc. Majority of such results use the incompressible Navier-Stokes system to model the fluid. We mention \cite{time,3dmesh,grandmont3,ruzicka,muhasch,BorSunNavierSlip,BorSun,trwa} for the existence of weak solutions, and \cite{strongzero,grandmont1,grandmont2,strong} for the existence of strong solutions. On the other hand, for the compressible Navier-Stokes system, there are only a handful of results which appeared in the last few years. First such result was due to Breit and Schwarzacher in \cite{compressible}, where they proved the existence of a weak solution for a shell interacting with a compressible viscous fluid. Later, Mitra proved the existence of a unique local regular solution in 2D, where the beam is viscoelastic in \cite{sourav}. Maity and Takahashi obtained a strong solution in $L^p-L^q$ framework for the interaction of heat-conducting fluid with a viscoelastic plate, and later Maity, Roy and Takahashi obtained a unique local regular solution in 3D \cite{comstrong}, where the structure is modeled by a wave equation and doesn't have any viscosity. Trifunovi\'{c} and Wang obtained the existence of a weak solution in 3D in \cite{trwa3}, where the structure is modeled by a nonlinear thermoelastic (or just elastic) plate, by constructing a time-continuous splitting approximation scheme that decouples the fluid and the structure.

 While weak solutions are in general global, or in the case of fluid-structure interaction exist up to a collision or self-collision of elastic structures, they are very often not unique and have a low regularity that comes from the energy of the problem. On the other hand, strong (or more regular) solutions are usually unique and posses higher regularity, but are often local in time. Thus, a natural question arises -  when do weak and strong solutions coincide if eminating from the same initial data, or in other words, when do weak solutions enjoy so called weak-strong uniqueness property? Note that this only gives a partial answer to the uniqueness problem, as it usually doesn't hold globally in time and requires a more regular initial data in order to ensure the existence of a strong solution in the first place. Some systems, like the incompressible Navier-Stokes equations, actually enjoy a better type of uniqueness properties. In particular, finite energy weak solutions in 2D are unique, while in 3D, a different type of a weaker uniquneness, so called Prodi-Serrin, is proven to hold. This means that a weak solution that has the integrability $\bb{u}\in L_t^a L_x^b$, for $a\in (2,\infty), b\in (3,\infty)$ such that $\frac{2}{a}+\frac{3}{b}=1$, is unique in the class of weak solutions. On the other hand, the compressible Navier-Stokes system is not as fortunate, as only the results concerning the weak-strong uniqueness are available so far (see \cite{comWSU2,comWSU1,comWSU0}), at least without any further simplifications.

The uniqueness of fluid-structure interaction problems has been studied mostly for the interaction of rigid bodies and incompressible viscous fluids. The following results are mentioned: the uniqueness in 2D case with no-slip condition \cite{uniqueness}, the uniqueness and the energy identity in 2D with slip condition \cite{bravin}, the weak-strong uniqueness with slip condition \cite{WSUrigid} and the Serri-Prodin uniqueness with slip or no-slip \cite{prodi-serrin}. The weak-strong uniqueness for compressible viscous fluid on moving domains, surrounding a rigid body and filling a rigid body were studied in \cite{movingdomains,comWSUrigid,comWSUin}, respectively. In the mentioned results, when comparing two solutions, the usual approach is to transform coordinates of both solutions to be expressed on a common fixed domain and then compare them there. One should note that to be even able to do this in the first place, the rigid bodies  corresponding to two different solutions need to be a priori isomorphic, while the fluid cavities need to be the same. Consequently, if one of the two solutions is more regular, then the rigid body and the fluid cavity for both of the solutions will need to have the same (higher) regularity. In other words, the  regularities of the boundaries are a priori fixed and do not play a role in the analysis. However, when one deals with the interaction of fluids and elastic plates (or shells), the regularity of fluid domains directly depends on the regularity of the plate displacements, and as such can be different for two different solutions. Uniqueness of such problems was first studied by Guidoboni, Guidorzi, and Padula in \cite{continuous}, where a 1D beam interacts with a 2D incompressible fluid. For this problem, they proved the continuous time dependence of solutions on initial data, and consequently the uniqueness of such weak solutions. Very recently, Schwarzacher and Sroczinski in \cite{WSinc} studied the interaction of incompressible viscous fluids and elastic plates, in 2D or 3D. Denoting the fluid velocities as $\bb{u}_i$ and the plate displacements as $w_i$, they proved that two finite energy weak solutions $(\bb{u}_1,w_1)$ and $(\bb{u}_2,w_2)$ coincide if $\partial_t \bb{u}_2 \in L_t^2 W_x^{-1,2}$ in 2D, or $\bb{u}_2 \in L_t^p W_x^{1,q}$ and $\partial_t \bb{u}_2 \in L_t^2 W_x^{-1,r}$ in 3D, for any $p>2,q>3$. The authors of this paper argue that the condition on $\partial_t \bb{u}_2$ was forgotten in \cite{continuous} and it needs to be additionally imposed. This is because when $(\bb{u}_2,w_2)$ is compared with $(\bb{u}_1,w_1)$, it is mapped onto a domain corresponding to the latter and this mapping (constructed in \cite[section 6]{continuous} and also used in \cite{WSinc}) is chosen to preserve the divergence-free condition in order to avoid the appearance of the pressure, but it doesn't ensure that $\partial_t \bb{u}_2$ is in the right dual space of divergence-free functions corresponding to the domain of the solution $(\bb{u}_1,w_1)$, so this is why this condition needs to be imposed additionally.

In this paper, we study the finite energy weak solutions (FEWS) to the interaction problem of compressible viscous fluids and elastic or viscoelastic plates. The notions of relative entropy and relative energy inequality are introduced. First, by utilizing the dissipation effects of the fluid onto the structure, we prove that the additional regularity $\Delta w \in L_t^2 H_x^s$ can be obtained, for $s = s(\gamma,d,\alpha)>0$, where $w$ is the displacement of the plate, $\gamma$ is the adiabatic constant, $d=2,3$ is the dimension of the fluid and $\alpha$ is the viscoelasticity constant.  This is done by following the ideas from \cite{muhasch}, where similar (but stronger) estimates were obtained for the interaction of a incompressible viscous fluid and a nonlinear Koiter shell. Afterwards, it is proved that all FEWS satisfy the relative energy inequality. Finally, we prove that all FEWS enjoy the weak-strong uniqueness property. The proof uses the strategy from \cite{comWSU2} (see also \cite{comWSU1,comWSU0}), and relies on using the relative entropy, which measures the distance between two solutions, and the relative energy inequality, which if tested by sufficiently regular solution of the same problem, gives the desired result by means of the Gronwall lemma. Since the fluid domains of two solutions are not the same in general, the coordinates of both solutions are transformed to the ones of a common fixed domain. The time and space derivatives are then replaced by the corresponding transformed derivatives which depend on the domain transformation, and thus on the displacement of the plate. This drastically increases the difficulty of the analysis. To close these estimates, additional integrability  is needed - $\nabla w \in L_t^2 L_x^\infty$, which is always ensured by the mentioned regularity $\Delta w \in L_t^2 H_x^s$ and the imbedding of Sobolev spaces (for any $s>0$). In the last part of the paper, the obtained results are extended to fluid-structure interaction models for which the structure is more nonlinear. More precisely, we first study the case where the structure is modeled by a semilinear plate model which generalizes Kirchhoff, von Karman and Berger plates (see Appendix for more details). Then, this equation is then replaced by a semilinear thermoelastic system and finally by a quasilinear thermoelastic system, thus furthering the frontiers of the approach presented here.

This paper is organized as follows. In section $\ref{sec2}$, the model we will study is introduced, then the notions of weak and strong solution for this problem are given, the relative entropy and relative energy inequality are introduced, and the main results are presented. In sections $\ref{proofofmain1}$ and $\ref{proofofmain2}$, the main results are proved. Finally in section 5, the main results are extended to interaction problems including certain nonlinear elastic/thermoelastic plates.

\section{The model}\label{sec2}

Here we deal with a compressible viscous fluid interacting with an elastic or viscoelastic plate, where the dimensions of the fluid and the plate are $d$ and $d-1$, respectively, for $d=2,3$. The vertical plate displacement is described by a scalar function $w: \Gamma\to \mathbb{R}$, where $\Gamma \subset \mathbb{R}$ is a connected bounded domain with a Lipschitz boundary. The fluid fills the domain defined as
\begin{eqnarray*}
\Omega^w(t) := \{ (X,z) : X \in \Gamma, -1< z < w(t,X)\}, \quad t\in[0,T].
\end{eqnarray*}
We will denote the graph of $w$ as $\Gamma^w (t) = \{ (X,z) : X \in \Gamma, z=w(t, X) \}$, while the entire rigid part of the boundary $\partial \Omega^w(t)$ will be denoted as $\Sigma := (\Gamma \times \{-1\}) \cup \partial\Gamma\times \{-1,0\}$. Finally, by a slight abuse of notation, for a time interval $I$ and a time dependent set $B(t)$, we will write
\begin{eqnarray*}%%%%%%%%%%%%%
	I\times B(t):= \bigcup_{t\in I}\{t\}\times B(t),
\end{eqnarray*}%%----------------------------%%
and the time-space cyllinders corresponding to our system will be denoted as
\begin{eqnarray*}%%%%%%%%%%%%%
		Q_T^w:=(0,T)\times \Omega^w(t), \quad\Gamma_T^w:=(0,T)\times\Gamma^w(t),  \quad \Gamma_T:=(0,T)\times\Gamma.
\end{eqnarray*}%%----------------------------%%
The governing equations for our coupled fluid-structure interaction problem read:\\

\noindent
\textbf{The (visco)elastic plate equation} on $\Gamma_T$:
\begin{equation}\label{structureeqs}
\partial_t^2 w+\Delta^2 w - \alpha \partial_t \Delta w=-S^w \bb{f}_{fl}\cdot \mathbf{e_3}; 
\end{equation}
\textbf{The compressible Navier-Stokes equations} on $Q_T^w$:
\begin{eqnarray}
\partial_t (\rho\mathbf{u}) + \nabla \cdot (\rho\mathbf{u}\otimes \mathbf{u})& =& -\nabla p(\rho) +\mathbb{S}(\nabla \bb{u}), \label{momeq}\\
\partial_t \rho + \nabla \cdot (\rho \mathbf{u}) &=& 0; \label{conteq}
\end{eqnarray}
\textbf{The fluid-structure coupling (kinematic and dynamic, resp.)} on $\Gamma_T$:
\begin{eqnarray}%%%%%%%%%%%%%
\partial_t w(t,X) \mathbf{e_3}&=& \mathbf{u}(t,X,w(t,X)),\label{kinc}\\
\bb{f}_{fl}(t,X)&=&\big[(-p(\rho)I+\nabla \cdot \mathbb{S}(\nabla \bb{u}))\nu^w\big](t,X,w(t,X)); \label{dync}
\end{eqnarray}%%----------------------------%%
\textbf{The boundary conditions (clamped and no-slip, resp.)}:
\begin{eqnarray}
w(t, x)=\partial_\nu w(t, x)&=&0, \quad \text{ on } (0,T)\times \partial\Gamma,\label{boundaryconditions1}\\
\mathbf{u}&=&0, \quad \text{ on}~~ (0,T) \times \Sigma \label{boundaryconditions2};
\end{eqnarray}
\vskip 0.1in
\noindent
\textbf{The initial data:}
\begin{eqnarray}\label{initialdata}
\rho(0,\cdot) = \rho_0, \quad(\rho\mathbf{u})(0, \cdot)=(\rho\bb{u})_0,\quad w(0,\cdot) = w_0,\quad\partial_t w(0,\cdot) = v_0. \\ \nonumber
\end{eqnarray}

\noindent
Here, $\alpha\geq 0$ is the viscoelastic coefficient, $S^w(t,X)$ is the Jacobian of the transformation from the Eulerian to the Lagrangian coordinates of the plate
\begin{eqnarray*}
S^w(t,X)=\sqrt{1+ |\nabla w|^2},
\end{eqnarray*}
$\bb{e}_3 = (0,0,1)$ (or $\bb{e}_3 = (0,1)$ if $d=2$), $\nu^w$ is the unit outward normal vector on $\Gamma^w$, $p$ is the pressure given by $p(\rho) = \rho^\gamma$, where\footnote{The lower bounds for $\gamma$ are imposed to ensure the existence of a weak solution (see Remark $\ref{thatremark}(3)$ for more details).}
\begin{enumerate}
	\item[] $\gamma>1$, if $d=2$ and $\alpha \geq 0$,
	\item[] $\gamma > \frac{12}{7}$, if $d=3$, and $\alpha = 0$,
	\item[] $\gamma>\frac{3}{2}$, if $d=3$ and $\alpha >0$,
\end{enumerate}
$\mathbb{S}(\nabla \bb{u})$ is the Newton stress tensor given by
\begin{eqnarray*}%%%%%%%%%%%%%
\mathbb{S}(\nabla \bb{u}):=\mu\nabla \bb{u}+(\mu+\lambda)(\nabla\cdot \bb{u})I,
\end{eqnarray*}%%----------------------------%%
with $\mu>0$ and $\lambda + \frac{2}{3}\mu> 0$, and $\nu$ is the normal vector on $\partial\Gamma$.

\subsection{Weak solution}
First, for any $1\leq p,q\leq \infty$, we define the following functional spaces for variable domains 
\begin{eqnarray*}%%%%%%%%%%%%%
    L^p(0,T; L^q(\Omega^w(t))):= &&\{ f\in L^1((0,T)\times \Omega^w(t)):\\ &&f(t,\cdot) \in L^q(\Omega^w(t)) \text{ for a.a }t\in(0,T),~ ||f(t,\cdot)||_{L^q(\Omega^w(t))}\in L^p(0,T)\}, 
\end{eqnarray*}%%----------------------------%%
and
\begin{eqnarray*}%%%%%%%%%%%%%
     && L^p(0,T; W^{1,q}(\Omega^w(t))):= \{ f\in  L^p(0,T; L^q(\Omega^w(t))): \nabla f \in  L^p(0,T; L^q(\Omega^w(t)))\},\\
    && L^p(0,T; W^{2,q}(\Omega^w(t))):= \{ f\in  L^p(0,T; W^{1,q}(\Omega^w(t))): \nabla^2 f \in  L^p(0,T; L^q(\Omega^w(t)))\}.
\end{eqnarray*}%%----------------------------%%
Now, to introduce the weak formulation to the problem $\eqref{structureeqs}-\eqref{initialdata}$, the following spaces that come from the energy inequality $\eqref{energy}$ are defined:
\begin{enumerate}
\item[] the structure displacement space\footnote{Here, $\alpha H^{1}(\Gamma_T):=\{f\in L^2(\Gamma_T): \alpha f\in H^1(\Gamma_T)\}$.}
\begin{eqnarray*}%%%%%%%%%%%%%
\mathcal{W}_S(0,T):=W^{1,\infty}(0,T;L^2(\Gamma)) \cap L^\infty(0,T;H_0^2(\Gamma))\cap \alpha H^1(\Gamma_T), 
\end{eqnarray*}%%----------------------------%%
\item[] the space for the fluid density
\begin{eqnarray*}%%%%%%%%%%%%%
\mathcal{W}_D(0,T):= L^\infty (0,T; L^\gamma(\Omega^w(t))).
\end{eqnarray*}%%----------------------------%%
\item[] the fluid velocity space
\begin{eqnarray*}%%%%%%%%%%%%%
\mathcal{W}_F(0,T):= L^\infty(0,T; L^2(\Omega^w(t)))\cap L^2(0,T; H^1(\Omega^w(t))),
\end{eqnarray*}%%----------------------------%%
\item[] and the coupled fluid-structure space
\begin{eqnarray*}%%%%%%%%%%%%%
\mathcal{W}_{F S}(0,T)= \{ (\bb{u}, w) \in \mathcal{W}_F(0,T) \times \mathcal{W}_S (0,T): \gamma_{|\Gamma^w(t)}\bb{u} = \partial_t w \bb{e}_3 \text{ and }\gamma_{|\Sigma}\bb{u}=0 \text{ for a.a. }t\in(0,T) \}.
\end{eqnarray*}%%----------------------------%%
\end{enumerate}
Here, $\gamma_{|\Gamma^w(t)}$ is the Lagrangian trace operator on $\Gamma^{w}(t)$, which is a continuous linear operator from $H^1(\Omega^w(t))$ to $H^s(\Gamma)$, for any $s<\frac{1}{2}$ (see \cite{Boris}). \\

The weak solution of the problem $\eqref{structureeqs}-\eqref{initialdata}$ is defined as:\footnote{See \cite{trwa3} for the derivation of this weak form for smooth solutions.}
\begin{mydef}\label{weaksolution}
We say that $(\rho ,\bb{u}, w) \in \mathcal{W}_D(0,T) \times \mathcal{W}_{FS}(0,T) $ is a \textbf{weak solution} of the problem $\eqref{structureeqs}$-$\eqref{initialdata}$ if:
\begin{enumerate}
\item The initial data $\rho_0, (\rho \bb{u})_0, w_0, v_0,\in L^\gamma(\Omega^{w_0})\times L^{\frac{2\gamma}{\gamma+1}}(\Omega^{w_0})\times H_0^2(\Gamma)\times L^2(\Gamma)$ and is assumed to satisfy the following compatibility conditions:
\begin{eqnarray*} %%%%%%%%%%%%%
\begin{aligned}
\rho_0>0,& \text{ on } \{(\rho\bb{u})_0 \neq 0\}, \\
\frac{(\rho \bb{u})_0^2}{\rho_0}&\in L^1(\Omega^{w_0}),\\
\partial_\nu w_0=w_0 =0, &\text{ on } \partial \Gamma, \\
w_0>-1, &\text{ on } \Gamma.
\end{aligned}
\end{eqnarray*}%%----------------------------%%
\item $\rho \in C_w (0,T; L^\gamma(\Omega^w(t)))$ and the continuity equation\footnote{Here, the density $\rho$ is extended by $0$ to $\mathbb{R}^3$.}
\begin{align}%%%%%%%%%%%%%
\int_0^T \int_{\mathbb{R}^3} \big( \rho \partial_t \varphi + \rho\bb{u}\cdot \nabla \varphi \big) = \int_{\mathbb{R}^3} \rho(T) \varphi(T) - \int_{\mathbb{R}^3} \rho_0 \varphi(0), \label{cont}
\end{align}%%----------------------------%%
holds for all $\varphi \in C^\infty([0,T]\times \mathbb{R}^3)$.
\item $\rho\bb{u} \in C_w (0,T; L^{\frac{2\gamma}{\gamma+1}}(\Omega^w(t)))$ and the coupled momentum equation
\begin{eqnarray}%%%%%%%%%%%%%
&&\int_{Q_T^w} \rho \bb{u} \cdot\partial_t \bb{q} + \int_{Q_T^w}(\rho \bb{u} \otimes \bb{u}):\nabla\bb{q} +\int_{Q_T^w} \rho^\gamma (\nabla \cdot \bb{q})- \int_{Q_T^w} \mathbb{S}(\nabla \bb{u}): \nabla \bb{q}\nonumber\\
&& +\int_{\Gamma_T} \partial_t w \partial_t \psi- \int_{\Gamma_T}\Delta w \Delta \psi- \alpha \int_{\Gamma_T} \partial_t \nabla w \cdot \nabla \psi\nonumber\\
&&=\int_{\Omega^w(T)}\rho(T)\mathbf{u}(T)\cdot\bb{q}(T)-\int_{\Omega^{w_0}}(\rho\mathbf{u})_0\cdot\bb{q}(0) + \int_\Gamma \partial_t w(T) \psi(T)- \int_\Gamma v_0 \psi(0), \quad \quad \quad \quad \quad \label{weaksolmom}
\end{eqnarray}%%----------------------------%%
holds for all $\bb{q} \in C^\infty([0,T]\times \overline{\Omega^w(t)})$ and $\psi\in C_0^\infty([0,T]\times\Gamma)$ such that $\bb{q}_{|\Sigma} =0$ and $\bb{q}(t,X,w(t,X))= \psi(t,X) \bb{e}_3$ on $\Gamma_T$.
\end{enumerate}
If, in addition $(\rho ,\bb{u}, w)$ satisfies the following energy inequality
\begin{eqnarray}%%%%%%%%%%%%%
&&\frac{1}{2} \bint_{\Omega^w(t)} (\rho|\bb{u}|^2) (t)+\frac{1}{\gamma-1} \bint_{\Omega^w(t)} \rho^\gamma(t) + \int_0^t\bint_{\Omega^w(t)} \mathbb{S}(\nabla \bb{u}): \nabla \bb{u} \nonumber \\
&&\quad+ \frac{1}{2}\int_\Gamma |\partial_t w(t)|^2+ \frac{1}{2} \int_\Gamma |\Delta w(t)|^2 +\alpha \int_0^t\bint_{\Gamma} |\partial_t \nabla w|^2 \nonumber \\
&& \leq \frac{1}{2} \bint_{\Omega^{w_0}} \frac{(\rho\bb{u})_0^2}{\rho_0 }+\frac{1}{\gamma-1} \bint_{\Omega^{w_0}} \rho_0^\gamma + \frac{1}{2}\int_\Gamma |v_0|^2+ \frac{1}{2} \int_\Gamma |\Delta w_0|^2, \label{energy}
\end{eqnarray}%%----------------------------%%
then we call $(\rho ,\bb{u}, w)$ a \textbf{finite energy weak solution}.
\end{mydef}

\begin{thm}(\cite[\textbf{Theorem 2.1}]{trwa3})
Let $\alpha = 0$ and $\gamma>\frac{12}{7}$ or $\alpha>0$ and $\gamma>\frac{3}{2}$. Then there exists a finite energy weak solution $(\rho,\bb{u},w)\in \mathcal{W}_D(0,T) \times \mathcal{W}_{FS}(0,T) $ in the sense of the Definition $\ref{weaksolution}$. The lifespan $T>0$ is either any $T<T^*$, where $T^*$ is the moment when the collision occurs $\min\limits_{X\in\Gamma}w(T^*,X)+1 = 0$, or $T=\infty$ if no collision occurs.
\end{thm}

\begin{rem}\label{thatremark}
	(1) In theorem above, the structure is modeled by a nonlinear thermoelastic plate model, but the same result holds for the linear elastic plate model given in $\eqref{structureeqs}$. In section 5, the nonlinear elastic/thermoelastic plate models are considered.\\
	(2) In \cite{compressible}, a weak solution in the sense very similar to the one in Definition $\ref{weaksolution}$ was constructed, for $\alpha = 0$ and $\gamma>\frac{12}{7}$. There however, the geometry is different - the elastic structure is a shell which deforms in the normal direction of the reference domain boundary.\\	
	(3) One of the main difficulties in studying the existence of such weak solutions is the weak $L^1$ convergence of the pressure. Since the structure regularity that comes from the energy is not strong enough to ensure that the fluid cavity is Lipschitz (and it also changes in time), the standard approach based on Bogovskii operator fails in this framework. The alternative here consists in proving the additional integrability of the pressure inside the domain in the standard way and that the $L^1$ norm of the pressure doesn't concentrate near the boundary, thus giving the equiintegrability sufficient to pass to the weak limit. The latter proof is very involved and it is based on constructing a test function which sort of represents the distance with respect to the boundary. As a consequence, this construction in \cite{compressible} and \cite{trwa3} varies greatly due to different geometries, as the latter has corners at $\partial \Gamma\times \{-1\}$ and $ \partial \Gamma\times \{0\}$. If $\alpha=0$, in order to close the estimates in the mentioned proof, one needs the condition $\gamma>12/7$ to compensate for the low integrability of the structure velocity $\partial_t w \in L_t^2 L_x^p$, for any $p<4$, which comes from the trace regularity. If however $\alpha>0$, then one can use the regularity $\partial_t w \in L_t^2 H_x^1$, so $\gamma>3/2$ suffices. This means that the same results as in \cite{compressible,trwa3} hold for $\alpha>0$ and $\gamma>3/2$, even though they weren't explicitly considered.\\

\end{rem} 
\subsection{Strong solution}
Here, we refer to known results concerning strong and regular solutions to the problem $\eqref{structureeqs}-\eqref{initialdata}$:
\begin{table}[h!]
	\begin{tabular}{|c|c|c|}
		\hline
		Reference, constants& Regularity of initial data & Regularity of solutions \\ 
		\hline
		\cite[\textbf{Theorem 1.1}]{sourav}& $\rho_0 \circ A_{w_0}\in H^2  $  &  $\rho \circ A_w \in C_t^0 H_x^2 \cap C_t^1 H_x^1 $ \\ 
	    $d=2$ &  $\bb{u}_0 \circ A_{w_0}\in H^3$ &  $\bb{u}\circ A_w \in L_t^2 H_x^3 \cap C_t^0 H_x^{\frac{5}{2}}\cap H_t^1 H_x^2 \cap C_t^1 H_x^1 \cap H_t^2 L_x^2$         \\		
		$\gamma\geq 1$&$w_0 = 0$   & $w\in C_t^0 H_x^{\frac{9}{2}} \cap H_t^1 H_x^4 \cap C_t^1 H_x^3 \cap H_t^2 H_x^2 $  \\
    	$\alpha>0$  &$v_0 \in H^3$ &	$\cap C_t^2 H_x^1 \cap H_t^3 L_x^2 $\\
		\hline
  		\cite[\textbf{Theorem 1.1}]{comstrong}& $	\rho_0 \circ A_{w_0} \in H^3$    &  $\rho \circ A_w \in H_t^1 H_x^3\cap W_t^{1,\infty} H_x^2$ \\ 
		$d=2,3$  & $\bb{u}_0 \circ A_{w_0} \in H^3$ &  $\bb{u} \circ A_w \in L_t^2 H_x^4  \cap H_t^1 H_x^2 \cap H_t^2 L_x^2$         \\		
		$\gamma\geq 1$ & $w_0 \in H^4$  & $w\in L_t^\infty H_x^4 \cap W_t^{1,\infty} H_x^3\cap  H_t^2 H_x^2   \cap  H_t^3 L_x^2$  \\
		$\alpha\geq 0$ & $v_0\in H^3$ &	\\
		\hline
		\cite[\textbf{Theorem 1.1}]{NSFFSI}& $	\rho_0 \circ A_{w_0} \in W^{1,q}$  &  $\rho \circ A_w \in W_t^p L_x^q\cap  L_t^p W_x^q$ \\ 
		$d=2,3$  &   $ \bb{u}_0 \circ A_{w_0} \in B_{q,p}^{2(2-1/p)}$ &  $\bb{u}\circ A_w \in L_t^p W_x^{2,q} \cap W_t^{1,p} L_x^q $     \\		
		$\gamma\geq 1$ & $w_0 \in B_{q,p}^{2(2-1/p)}$  & $w\in L_t^p W_x^{4,q}\cap W_t^{2,p}L_x^q $  \\
		$\alpha>0$ & $ v_0\in  B_{q,p}^{2(2-1/p)}$ & $p\in(1,\infty), q\in(d,\infty), \frac{1}{p}+\frac{1}{2q}\neq \frac{1}{4-d}$	\\
		\hline
		
	\end{tabular}
\caption{List of strong and regular solutions to $\eqref{structureeqs}-\eqref{initialdata}$. The regularity of $\rho,\bb{u}$ is given in ALE coordinates, by means of the mapping $A_w$ from $\eqref{ALE}$, and the compatibility conditions for the initial data are omitted as they are quite natural and can be found in the respective references.  
}
\end{table}

\begin{rem}
	(1) In the above table, $B_{q,p}^{2(2-1/p)}$ stands for Besov spaces (see \cite{adams}). \\
	(2) In \cite{comstrong,sourav}, the plate domain $\Gamma$ is periodic, i.e. $\Gamma= (R/ L_1\mathbb{Z})\times (R/ L_2\mathbb{Z})$ and $\Gamma= (R/ L_1\mathbb{Z})$ , respectively, for some $L_1,L_2>0$, and  $\rho,\bb{u},w$ are periodic in horizontal coordinates as well. However, in \cite{NSFFSI}, the reference fluid domain is smooth and contains a flat part which is deformable and constitutes a viscoelastic plate.\\
	(3) In \cite{NSFFSI}, the fluid conducts heat and it is governed by the Navier-Stokes-Fourier system, but the same result holds for the barotropic case used in this paper, as it was pointed out in \cite[Remark 1.3(5)]{NSFFSI}.\\
	(4) In \cite{comstrong}, the structure is actually governed by a wave equation. Nevertheless, the same result can be obtained for the plate equation, both for $\alpha= 0$ and $\alpha>0$, as these cases  have more regularity. However, in \cite{NSFFSI,sourav}, $\alpha>0$ and the same results cannot be obtained in the same way if $\alpha=0$. This is because viscoelasticity ensures maximal regularity \cite{NSFFSI} and the analiticity \cite{sourav} properties, which play a key role in the respective approaches.
\end{rem}

\subsection{Relative entropy and relative energy inequality}
Relative energy (or entropy) inequality in the context of compressible Navier-Stokes system was first introduced by Germain in \cite{comWSU0}, and it was used as a tool in proving the weak-strong uniqueness property for a sub-class of finite energy weak solutions which have have a $L_t^p L_x^q$-integrable gradient of density (for certain $p,q>1$). Note that the latter condition is quite unnatural and existence of weak solutions within this class is non-existant in the literature. Later, Feireisl, Novotn\'{y} and Sun in \cite{comWSU1} introduced suitable weak solutions - a sub-class of finite energy weak solutions that satisfy the relative energy inequality, then proved their existence and showed the weak-strong uniqueness within this class, thus removing the additional regularity on density imposed in \cite{comWSU0}. Finally, Feireisl, Novotny and Jin in \cite{comWSU2} introduced the relative entropy in the class of finite energy weak solutions and gave a definitive answer to the weak-strong uniqueness question, proving that all finite energy weak solutions are in fact suitable weak solutions and thus have the weak-strong uniqueness property. Here, we extend the notions of relative entropy and relative energy inequality to our fluid-structure interaction model.\\

Let $(\rho_1,\bb{u}_1,w_1)$ be a finite energy weak solution in the sense of Definition $\ref{weaksolution}$. For any $\rho_2,\bb{u}_2 \in C^\infty([0,T]\times \overline{\Omega^{w_1}(t)})$ and $w_2 \in C_0^\infty([0,T]\times\Gamma)$ such that ${\bb{u}_2}_{|\Sigma} =0$ and $\bb{u}_2(t,X,w_2(t,X))= \partial_t w_2(t,X) \bb{e}_3$ on $\Gamma_T$, we define the \textbf{relative entropy} with respect to $(\rho_2,\bb{u}_2,w_2)$ as \\
\begin{eqnarray}%%%%%%%%%%%%%
\mathcal{E}\Big((\rho_1,\bb{u}_1,w_1)\Big|(\rho_2,\bb{u}_2,w_2)\Big)(t)  &=& \frac{1}{2} \int_{\Omega^{w_1}(t)} (\rho_1 |\bb{u}_1 - \bb{u}_2|^2)(t) + \frac{1}{\gamma-1}\int_{\Omega^{w_1}(t)} \big( \rho_1^\gamma - \gamma \rho_2^{\gamma-1}(\rho_1 - \rho_2) - \rho_2^{\gamma} \big)(t) \nonumber \\[2mm]
&&\quad + \frac{1}{2}\int_\Gamma |\partial_t w_1 - \partial_t w_2|^2(t) + \frac{1}{2}\int_\Gamma |\Delta w_1 - \Delta w_2|^2(t),\nonumber \\ \nonumber
\end{eqnarray}%%----------------------------%%
and the \textbf{relative energy inequality} with respect to $(\rho_2,\bb{u}_2,w_2)$ as\\
\begin{eqnarray}%%%%%%%%%%%%%
&&\mathcal{E}\Big((\rho_1,\bb{u}_1,w_1)\Big|(\rho_2,\bb{u}_2,w_2)\Big)(t)  + \int_{0}^t\int_{\Omega^{w_1}(t)}  \mathbb{S} (\nabla \bb{u}_1 - \nabla \bb{u}_2):(\nabla\bb{u}_1 - \nabla \bb{u}_2) + \alpha\int_0^t\int_{\Gamma}  |\partial_t \nabla w_1 - \partial_t \nabla w_2|^2 \nonumber\\
&&\leq \mathcal{E}\Big((\rho_1,\bb{u}_1,w_1)\Big|(\rho_2,\bb{u}_2,w_2)\Big)(0) +\int_0^t\mathcal{R}\Big( \rho_1,\bb{u}_1,w_1, \rho_2,\bb{u}_2,w_2\Big),\label{relen} \\ \nonumber
\end{eqnarray}%%----------------------------%%
where the \textbf{remainder term} reads\\
\begin{align}%%%%%%%%%%%%%
&\mathcal{R}\Big( \rho_1,\bb{u}_1,w_1, \rho_2,\bb{u}_2,w_2\Big) :=  \int_{\Omega^{w_1}(t)}  \mathbb{S}( \nabla \bb{u}_2):(\nabla\bb{u}_2 - \nabla \bb{u}_1)+\int_{\Omega^{w_1}(t)} \rho_1 (\partial_t \bb{u}_2 + \bb{u}_1\cdot \nabla \bb{u}_2)\cdot(\bb{u}_2 - \bb{u}_1) \nonumber\\
&\quad+ \frac{\gamma}{\gamma-1} \int_{\Omega^{w_1}(t)}  \Big[(\rho_2\bb{u}_2 - \rho_1 \bb{u}_1)\cdot \nabla (\rho_2^{\gamma-1}) + (\rho_1 - \rho_2)\partial_t(\rho_2^{\gamma-1})\Big]+ \int_{\Omega^{w_1}(t)} (\rho_2^\gamma - \rho_1^\gamma)(\nabla \cdot \bb{u}_2) \nonumber\\
&\quad+ \int_{\Gamma^{w_1}(t)}(\partial_t w_2 - \partial_t w_1) \rho_2^\gamma \nu^{w^1} \cdot \bb{e}_3-\int_{\Gamma} (\partial_t w_2 - \partial_t w_1) \partial_t^2 w_2 -\int_{\Gamma} (\partial_t w_2 - \partial_t w_1) \Delta^2 w_2\nonumber\\
&\quad+\alpha\int_{\Gamma} \partial_t \Delta w_2 (\partial_t  w_1 - \partial_t w_2) . \nonumber
\end{align}%%----------------------------%%
It is straightforward to see that by choosing $\rho_2 = 0, \bb{u}_2 = 0, w_2 = 0$, the relative energy inequality reduces to the energy inequality $\eqref{energy}$. Also, if one chooses $w_1 = w_2 = 0$, then relative entropy and relative energy inequality take the form of the ones defined in \cite{comWSU2}. \\

\subsection{Main results}
First, the following characterization of finite energy weak solutions is given:
\begin{thm}(\textbf{Main result I})\label{main2}
	Let $(\rho,\bb{u},w)$ be a \textbf{finite energy weak solution} in the sense of Definition $\ref{weaksolution}$ defined on $(0,T)$ such that $\min_{(0,T)\times \Omega^w(t)}1+w\geq c>0$. Then, one has the following:
	\begin{enumerate}
		\item[(1)] $\Delta w \in L^2(0,T; H^s(\Gamma))$, where
		\begin{enumerate}
			\item[2D:] $\alpha = 0$, $\gamma>1$ and $0<s<\min\big\{ \frac{1}{2} - \frac{1}{2\gamma}, \frac{1}{4} \big\}$;\\ $\alpha > 0$, $\gamma>1$ and $0<s<\min\big\{ \frac{3}{4} - \frac{1}{2\gamma}, \frac{1}{2} \big\}$;
			\item[3D:] $\alpha = 0$, $\gamma>\frac{12}{7}$ and $0<s<\min \big\{ \frac{7}{12} - \frac{1}{\gamma},\frac{1}{4} \big\}$;\\
			$\alpha >0$: $\gamma>\frac{3}{2}$ and $0<s<\min \big\{ \frac{2}{3} - \frac{1}{\gamma}, \frac{1}{2} \big\}$.
		\end{enumerate}
Consequently, $\nabla w \in L^2(0,T; L^\infty(\Gamma))$ in all the above cases;
\item[(2)] $(\rho,\bb{u},w)$ satisfies the relative energy inequality $\eqref{relen}$.
\end{enumerate}
\end{thm}

The second result of this paper states that all finite energy weak solutions enjoy the weak-strong uniqueness property in the following sense:

\begin{thm}(\textbf{Main result II})\label{main1}
Let $(\rho_1,\bb{u}_1,w_1)$ be a \textbf{finite energy weak solution} in the sense of Definition $\ref{weaksolution}$ defined on $(0,T_1)$ such that
\begin{eqnarray*}%%%%%%%%%%%%%
    \min\limits_{(0,T_1)\times \Omega^{w_1}(t)} 1+w_1 \geq c_1>0
\end{eqnarray*}%%----------------------------%%
and let $(\rho_2,\bb{u}_2,w_2)$ be a \textbf{strong solution} of the problem $\eqref{structureeqs}$-$\eqref{initialdata}$ such that
\begin{eqnarray*}%%%%%%%%%%%%%
	&&\rho_2 \in L^2(0,T_2; W^{1,q}(\Omega^{w_2}(t))) \cap H^1(0,T_2; L^q(\Omega^{w_2}(t))), \\
	&&\bb{u}_2 \in L^2(0,T_2; W^{2,q}(\Omega^{w_2}(t))) \cap H^1(0,T_2; L^q(\Omega^{w_2}(t))),\\
	&&w_2 \in L^2(0,T_2; H^4(\Gamma)) \cap H^2(0,T_2; L^2(\Gamma)),
\end{eqnarray*}%%----------------------------%%
and
\begin{eqnarray*}%%%%%%%%%%%%%
	0< \inf\limits_{(0,T_2)\times \Omega^{w_2}(t)} \rho_2\leq \sup\limits_{(0,T_2)\times \Omega^{w_2}(t)} \rho_2 <\infty, \qquad \min\limits_{(0,T_2)\times \Omega^{w_2}(t)} 1+w_2 \geq c_2>0.
\end{eqnarray*}%%----------------------------%%
If
\begin{eqnarray*}%%%%%%%%%%%%%
\rho_1(0,\cdot) = \rho_2(0,\cdot), \quad \bb{u}_1(0,\cdot)=\bb{u}_2(0,\cdot),\quad w_1(0,\cdot)= w_2(0,\cdot), \quad v_1(0,\cdot) =v_2(0,\cdot),
\end{eqnarray*}%%----------------------------%%
 and
\begin{enumerate}
\item[2D:] $\alpha \geq 0$, $\gamma>1$ and $q>2$;
\item[3D:] $\alpha = 0$, $\gamma\in \big(\frac{12}{7},2\big]$, $\nabla \bb{u}_2 \in L^\infty((0,T_2)\times \Omega^{w_2}(t))$ and $q>\frac{4\gamma}{3\gamma-4}$;\\[2mm]
$\alpha = 0$, $\gamma>2$ and $q>4$;\\[2mm]
$\alpha >0$, $\gamma>\frac{3}{2}$ and $q>\max\{3, \frac{6\gamma}{5\gamma-6} \}$,
\end{enumerate}
then
\begin{eqnarray*}%%%%%%%%%%%%%
(\rho_1,\bb{u}_1, w_1 ) \equiv (\rho_2,\bb{u}_2,w_2), \quad \text{a.e. in } (0,T_2)\times \Omega^{w_2}(t).
\end{eqnarray*}%%----------------------------%%
\end{thm}
${}$\\

The following extension of main results is obtained in section 5:
\begin{cor}
Theorems $\ref{main1}$ and $\ref{main2}$ hold in  the case where the linear plate model $\eqref{structureeqs}$ is replaced by certain nonlinear plate models.
\end{cor}

\begin{rem}
The purpose of the regularity estimate in Theorem $\ref{main2}(1)$ is threefold. First, it gives us additional properties for finite-energy weak solutions which can be used on their own in further applications. Second, it ensures that the domain $\Omega^w(t)$ is always Lipschitz (with values in $L_t^2$), which allows us to use the divergence theorem and Raynolds transport theorem in the proof of Theorem $\ref{main2}(2)$. Finally, it is used in the proof of Theorem $\ref{main1}$ to control 7th and 10th terms on the right-hand side of inequality $\eqref{estimate}$.
\end{rem}

\section{Proof of Theorem $\ref{main2}$}\label{proofofmain1}
The following proof will be carried in two parts.

\subsection{Part 1: the regularity estimates - Theorem $\ref{main2}(1)$}
We start with introducing fractional Sobolev and Nikolskii functional spaces, respectively (see \cite{adams} and \cite{triebel} for more details):
\begin{mydef}
	For $\alpha\in(0,1)$, $q\in (1,\infty)$, we say that $g \in \mathcal{W}^{\alpha,q}(\Gamma)$ if its norm satisfies
	\begin{eqnarray*}%%%%%%%%%%%%%
		|| g||_{\mathcal{W}^{\alpha,q}(\Gamma)}^q:= \Big( \int_\Gamma \int_\Gamma \frac{|g(x)-g(y)|^q}{|x-y|^{n+\alpha q}} dxdy \Big)^{\frac{1}{q}} + \Big( \int_\Gamma |g(x)|^q dx \Big)^{\frac{1}{q}} < \infty,
	\end{eqnarray*}%%----------------------------%%
and that $g \in \mathcal{N}^{\alpha,q}(\Gamma)$ if its norm satisfies
	\begin{eqnarray*}%%%%%%%%%%%%%
		|| g||_{\mathcal{N}^{\alpha,q}(\Gamma)}^q:= \sup\limits_{i\in \{1,...,d\}} \sup\limits_{h>0} \Big( \int_{\Gamma_h} \Big| \frac{g(x+h \bb{e}_i) - g(x)}{h^\alpha} \Big|^q dx \Big)^{\frac{1}{q}} + \Big( \int_\Gamma |g(x)|^q dx \Big)^{\frac{1}{q}}< \infty,
	\end{eqnarray*}%%----------------------------%%
	where $\bb{e}_i$ its the $i$-th unit vector and $\Gamma_h = \{x \in \Gamma: \text{dist}(x,\partial \Gamma)>h \}$. Given any $0<\alpha<\beta<1$ and $g \in \mathcal{N}^{\beta,q}(\Gamma)$, the following imbedding is true
	\begin{eqnarray*}%%%%%%%%%%%%%
		||g||_{L^{\frac{nq}{n-\alpha q}}(\Gamma)} \leq C_1||g||_{\mathcal{W}^{\alpha,q}(\Gamma)} \leq C_2||g||_{\mathcal{N}^{\beta,q}(\Gamma)}\leq C_3 ||g||_{\mathcal{W}^{\beta,q}(\Gamma)} , \quad \text{ for  }\alpha q<n ,
	\end{eqnarray*}%%----------------------------%%
	or
	\begin{eqnarray*}%%%%%%%%%%%%%
		||g||_{C^{\alpha-\frac{n}{q}}(\Gamma)}\leq C_1||g||_{\mathcal{W}^{\alpha,q}(\Gamma)} \leq C_2 ||g||_{\mathcal{N}^{\beta,q}(\Gamma)} \leq C_3||g||_{\mathcal{W}^{\alpha,q}(\Gamma)},\quad \text{ for  }\alpha q>n.
	\end{eqnarray*}%%----------------------------%%
	Finally, for $s,h>0$ and a function $q : \Gamma \to \mathbb{R}$, we introduce the quotient
	\begin{eqnarray*}%%%%%%%%%%%%%
		D_{h,\bb{e}}^s(q)(y):= \frac{q(y+h \bb{e}) - q(y)}{h^s}, ~~\text{for any (unit) vector } \bb{e} \in \mathbb{R}^2.
	\end{eqnarray*}%%----------------------------%%
	Since the direction of the unit vector will not be relevant in most cases, we will usually omit it and just write $D_h^s(q)(y):=D_{h,\bb{e}}^s(q)(y)$ instead, for any unit vector $\bb{e}$. 
	
\end{mydef}
The following notation will be useful in the upcoming analysis:

\begin{mydef}
	For a given $1<b \leq \infty$ and a domain $E$, denote by
	\begin{eqnarray*}%%%%%%%%%%%%%
		&L^{b^-}(E) := \bigcap_{p<b} L^p(E), \quad \quad W^{{a^-},b}(E) := \bigcap_{s<a} W^{s,p}(E),& \\[2mm]
		&W^{a,b^-}(E) := \bigcap_{p<b} W^{s,p}(E), \quad \quad
		W^{a^-,b^-}(E) := \bigcap_{s<a,p<b} W^{s,p}(E).&
	\end{eqnarray*}%%----------------------------%%
	We will write for any $1< b \leq \infty$,
	\begin{eqnarray*}%%%%%%%%%%%%%
		|| f ||_{L^{b^-}(E)}\leq D,
	\end{eqnarray*}%%----------------------------%%
	for a constant $D\geq 0$, if for all $1\leq b' <b$, there exists a constant $C(b')$ such that
	\begin{eqnarray*}%%%%%%%%%%%%%
		|| f ||_{L^{b'}(E)}\leq C(b')D.
	\end{eqnarray*}%%----------------------------%%
	and similarly for the Sobolev spaces $W^{a^-,b}(E)$, $W^{a,b^-}(E)$ and $W^{a^-,b^-}(E)$.
\end{mydef}
\begin{rem}
(1) Note that $|| f ||_{L^{b^-}(E)}$ is just a notation, not a proper norm;\\
(2) Wherever we have $|| f ||_{L^{b^-}(E)}$ appearing in an estimate, it can be replaced with $|| f ||_{L^{p}(E)}$ followed by "for $p$ close enough to $b$". By using this notation, we avoid writing this below the majority of estimates, especially in the case when there is more than one being used. Moreover, when $b^-$ appears in the estimates (instead of somewhere afterwards), it makes lengthy estimates easier to follow.
\end{rem}

\begin{rem}
	(1) The following proof is the adaptation to the compressible case of \cite[Theorem 1.2]{muhasch}, where similar estimates were obtained in the context of incompressible viscous fluid interacting with a quasilinear elastic Koiter shell;\\
	(2) We will only prove the estimate for the 3D/2D case, since the estimates in the 2D/1D case can be obtained more easily.
\end{rem}

We separately observe two cases for the domain $\Gamma$, first the periodic boundary domain and then a bounded Lipschitz domain.

\subsubsection{The periodic boundary case}
Here, we assume that $\Gamma= (R/ L_1\mathbb{Z})\times (R/ L_2\mathbb{Z})$ for $L_1,L_2>0$. Let $\bb{R}_w$ be an extension operator defined as
\begin{eqnarray*}%%%%%%%%%%%%%
	\bb{R}_w: f(t,X) \mapsto f(t,X)\frac{z+1}{w(t,X)+1}\bb{e}_3.
\end{eqnarray*}%%----------------------------%%

\noindent
Fix $h>0$ and in the coupled momentum equation $\eqref{weaksolmom}$, choose $(\bb{q},\psi) = (\bb{R}_w\big[D_{-h}^s D_h^{s}[w]\big], D_{-h}^s D_h^{s} [w])$ to obtain:
\begin{align}
&||D_h^s \Delta w||_{L^2(\Gamma_T)}^2 \nonumber \\
&= -\bint_{Q_T^w}(\rho \bb{u}\otimes \bb{u}):
\nabla \bb{R}_w \big[D_{-h}^s D_h^{s}[w]\big] -\bint_{Q_T^w} \rho\bb{u} \cdot
\partial_t \bb{R}_w\big[D_{-h}^s D_h^{s} [w] \big]
-\bint_{Q_T^w}\rho^\gamma \big(\nabla \cdot \bb{R}_w\big[D_{-h}^s D_h^{s}[w]\big]\big)\nonumber \\
&\quad+\mu\bint_{Q_T^w} \nabla \bb{u}:\nabla\bb{R}_w\big[D_{-h}^s D_h^{s}[w]\big]+(\mu+\lambda)\bint_{Q_T^w} (\nabla \cdot \bb{u})(\nabla \cdot \bb{R}_w\big[D_{-h}^s D_h^{s} [w]\big])\nonumber\\[2mm]
&\quad
+ ||D_{h}^s\partial_t w||_{L^2(\Gamma_T)} + \alpha \bint_{\Gamma_T} \partial_t  \nabla w \cdot \nabla D_{-h}^s D_h^{s}[w ] \nonumber\\[2mm]
&\quad +\int_{\Omega^w(T)} \rho(T) \mathbf{u}(T) \cdot \bb{R}_w\big[D_{-h}^s D_h^{s}[w(T)]\big]-\int_{\Omega^{w_0}} (\rho \mathbf{u})_0 \cdot \bb{R}_w\big[D_{-h}^s D_h^{s}[w_0]\big]) \nonumber \\
&\quad +\int_\Gamma \partial_t w(T) D_{-h}^s D_h^{s} w(T)-\int_\Gamma v_0 D_{-h}^s D_h^{s} w_0, & \label{aaaa}
\end{align}
where we used the following commutative property
\begin{eqnarray*}%%%%%%%%%%%%%
	-\int_{\Gamma} f D_{-h}^s D_h^s g = \int_{\Gamma}D_h^s f D_h^s g, \quad f,g \in L^2(\Gamma).
\end{eqnarray*}%%----------------------------%%
The proof will follow once we show that right-hand side of $\eqref{aaaa}$ can be bounded by a constant independent of $h$. We start with
\begin{eqnarray*}%%%%%%%%%%%%%
	&&\bint_{Q_T^w}(\rho \bb{u}\otimes \bb{u}):
	\nabla \bb{R}_w\big[D_{-h}^s D_h^{s}[w]\big] \\	
	&&\leq C|| \rho||_{L^\infty(0,T;L^\gamma(\Omega^w(t)))} ||\bb{u}||_{L^2(0,T; L^{6-}(\Omega^w(t)))}^2  \left|\left|\nabla \bb{R}_w\big[D_{-h}^s D_h^{s}  [w]\big]\right|\right|_{L^\infty(0,T; L^p(\Omega^w(t)))}\\[2mm]
	&&\leq C \Big|\Big| D_{-h}^s D_h^{s} [\nabla w] \frac{z+1}{w+1}  + D_{-h}^s D_h^{s} [w] \nabla \Big(\frac{z+1}{w+1}\Big) \Big|\Big|_{L^\infty(0,T; L^p(\Omega^w(t)))}\\[2mm]
	&&\leq C \left(|| D_{-h}^s D_h^{s} [\nabla w]||_{L^\infty(0,T; L^p(\Gamma))} +|| D_{-h}^s D_h^{s}  [w]||_{L^\infty(0,T; L^{p^-}(\Gamma))}\Big|\Big|\nabla \Big(\frac{z+1}{w+1}\Big) \Big|\Big|_{L^\infty(0,T; L^{\infty^-}(\Omega^w(t)))} \right) \\[2mm]
	&& \leq  C \big( ||\nabla w||_{L^\infty(0,T; W^{2s,p}(\Gamma))}+ || w||_{L^\infty(0,T; W^{2s,p^-}(\Gamma))} ||\nabla w ||_{L^\infty(0,T; L^{\infty^-}(\Gamma))}\big) \\[2mm]
	&&\leq C|| w||_{L^\infty(0,T; H^{2}(\Gamma))} \leq C,
\end{eqnarray*}%%----------------------------%%
where we have used the H\"{o}lder's inequality to obtain the second line which constraints $p>\frac{3\gamma}{2\gamma-3}$ and the imbedding $H^{1}(\Gamma)\hookrightarrow W^{2s,p}(\Gamma)$ to obtain the last line which constraints $s<\frac12$ and $\frac{4}{2-(1-2s)2}=\frac1s\geq p$, so
\begin{eqnarray*}%%%%%%%%%%%%%
   s< \frac12 ~\text{ and }~ \frac{1}{s}\geq p>\frac{3\gamma}{2\gamma-3} \implies 0<s<\min \left\{ \frac{2}{3}-\frac{1}{\gamma},\frac12 \right\}.
\end{eqnarray*}%%----------------------------%%
Next, 
\begin{eqnarray*}%%%%%%%%%%%%%
	&&\bint_{Q_T^w} \rho\bb{u} \cdot
	\partial_t \bb{R}_w\big[D_{-h}^s D_h^{s} [w]\big] 	\leq ||\rho\bb{u}||_{L^2(0,T; L^{(\frac{6\gamma}{5\gamma-6})^-}(\Omega^w(t)))} \left|\left| 	\partial_t \bb{R}_w\big[D_{-h}^s D_h^{s} [w]\big] \right|\right|_{L^2(0,T;L^p(\Omega^w(t)))} \\
	&&\leq C \Big|\Big| D_{-h}^s D_h^{s}  [\partial_t w] \frac{z+1}{w+1}- D_{-h}^s D_h^{s} [w] \frac{z+1}{(w+1)^2} \partial_t w    \Big|\Big|_{L^2(0,T;L^p(\Omega^w(t)))}\\[2mm]
	&&\leq C\big(|| D_{-h}^s D_h^{s}  [\partial_t w]     ||_{L^2(0,T;L^p(\Gamma))}+|| D_{-h}^s D_h^{s}  [ w] \partial_t w     ||_{L^2(0,T;L^p(\Gamma))} \big) \\
	&&\leq C\big(|| \partial_t w  ||_{L^2(0,T;W^{2s,p}(\Gamma))}+|| D_{-h}^s D_h^{s}[w] ||_{L^\infty(0,T;L^\infty(\Gamma))} || \partial_t w ||_{L^2(0,T;L^p(\Gamma))} \big) \\
	&&\leq C(1+ || D_{-h}^s D_h^{s}[w] ||_{L^\infty(0,T;W^{1,q}(\Gamma))})\\
	&&\leq  C(1+  ||w ||_{L^\infty(0,T;W^{1+2s,q}(\Gamma))}) \leq C(1+|| w||_{L^\infty(0,T; H^{2}(\Gamma))})\leq C,
\end{eqnarray*}%%----------------------------%%
where we have used the H\"{o}lder's inequality in the first inequality which constraints $p>\frac{6\gamma}{5\gamma-6}$, the imbedding in the line before the last one which constraints $q>2$ (but close to $2$), the imbedding in the last line which constraints  $s<\frac12$ so that $H^{2}(\Gamma)\hookrightarrow W^{1+2s,q}(\Gamma)$ and
\begin{enumerate}
    \item Case $\alpha=0$: in order to bound
    \begin{eqnarray*}%%%%%%%%%%%%%
        || \partial_t w  ||_{L^2(0,T;W^{2s,p}(\Gamma))} \leq C||\partial_t w||_{L^2(0,T;H^{(\frac12)^{-}}(\Gamma))}\leq C
\end{eqnarray*}%%----------------------------%%
    we need the following imbedding to hold $H^{(\frac12)^-}(\Gamma)\hookrightarrow W^{2s,p}(\Gamma)$ which constraints $s<\frac14$, $p\leq \frac{4}{2-(1/2-2s)2}=\frac{4}{1+4s}$ and $p<4$, so 
    \begin{eqnarray*}%%%%%%%%%%%%%
        s<\frac14 ~\text{ and }~ \frac{4}{1+4s} \geq p> \frac{6\gamma}{5\gamma-6} \implies 0<s< \min\Big\{\frac{7}{12}-\frac{1}{\gamma},\frac{1}{4}\Big\};
\end{eqnarray*}%%----------------------------%%
\item Case $\alpha>0$: in order to bound
\begin{eqnarray*}%%%%%%%%%%%%%
    || \partial_t w  ||_{L^2(0,T;W^{2s,p}(\Gamma))} \leq C||\partial_t w||_{L^2(0,T;H^{1}(\Gamma))}\leq C
\end{eqnarray*}%%----------------------------%%
we need $H^{1}(\Gamma)\hookrightarrow W^{2s,p}(\Gamma)$ to hold which constraints $s< \frac12$, $p\leq \frac{4}{2-(1-2s)2}=\frac{1}{s}$ and $p<\infty$, so
\begin{eqnarray*}
 s<\frac12 ~\text{ and }~  \frac{1}{s} \geq p> \frac{6\gamma}{5\gamma-6} \implies 0<s< \min\Big\{\frac{5}{6}-\frac{1}{\gamma},\frac{1}{2}\Big\}.
\end{eqnarray*}%%----------------------------%%
\end{enumerate}
Next, notice that we can explicitly calculate on $Q_T^w$
\begin{eqnarray*}%%%%%%%%%%%%%
	\nabla \cdot \bb{R}_w\big[D_{-h}^s D_h^{s}[w]\big] = D_{-h}^s D_h^{s}[w] \frac{1}{w+1},
\end{eqnarray*}%%----------------------------%%
so one can bound
\begin{eqnarray*}%%%%%%%%%%%%%
	&&\bint_{Q_T^w}\rho^\gamma \big(\nabla \cdot \bb{R}_w\big[D_{-h}^s D_h^{s}[w]\big]\big) \leq C|| \rho^\gamma||_{L^\infty(0,T;L^1(\Omega^w(t)))} \Big|\Big| D_{-h}^s D_h^{s}[w]\frac{1}{w+1}\Big|\Big|_{L^1(0,T; L^\infty(\Gamma))} \\[2mm]
	&&\leq  C || D_{-h}^s D_h^{s} [w] ||_{L^2(0,T; W^{1,p}(\Gamma))} \leq C || w||_{L^2(0,T; W^{1+2s,p}(\Gamma))} \leq C || w||_{L^2(0,T; H^{2}(\Gamma))} \leq C,
\end{eqnarray*}%%----------------------------%%
by the imbedding of Sobolev spaces $L^\infty(\Gamma)\hookrightarrow W^{1,p}(\Gamma)$ and $H^2(\Gamma)\hookrightarrow W^{2s,p}(\Gamma)$, which constraints $p>2$ (but arbitrarily close to $2$) and $s<\frac12$. The remaining terms can be bounded in a similar fashion as the previous ones, the proof is finished.

\subsubsection{Bounded domain case}
In this case, we choose $\Gamma$ to be a bounded and Lipschitz domain. There are three parts to this proof, first one are the interior estimates, second one are the boundary estimates for the flat boundary case and the third one are the boundary estimates for general boundary case. The detailed calculation is omitted here and only the main steps are presented. An interested reader is referred to \cite[Section 6.3]{evans}, where the same steps are taken to obtain the improved regularity estimates for the second order elliptic non-homogeneous equation.\\

\noindent
\textbf{Step 1 - interior estimates}. Let $V\Subset \Gamma$ and fix $h>0$ such that $h<\text{dist}(V,\partial\Gamma)$. Let $\phi\in C_c^\infty(\Gamma)$ be a cut-off function such that $0\leq\phi\leq1$ on $\Gamma$ and $\phi=1$ on $V$. We choose $(\bb{q},\psi) = (\bb{R}_w\big[D_{-h}^s D_h^{s}[w\phi]\big], D_{-h}^s D_h^{s} [w\phi])$ in the coupled momentum equation $\eqref{weaksolmom}$. We obtain an identity almost the same to the one in $\eqref{aaaa}$ with some additional terms involving the derivatives of $\phi$, so by estimating in the same way as in the periodic domain case, we obtain that $||\phi D_h^s \Delta w||_{L^2(\Gamma_T)}^2\leq C(\text{dist}(V,\partial\Gamma))$. Since $V$ was arbitrary, this gives us $\Delta w \in L^2(0,T;H_{loc}^s(\Gamma))$, where $s$ satisfies the conditions given in Theorem $\ref{main2}(1)$.\\

\noindent
\textbf{Step 2 - boundary estimates (flat boundary case)}. Here, we assume that that we are dealing with a flat part of the boundary. More precisely, let w.l.o.g. $(x_0,y_0)=X_0\in \partial\Gamma$ and
\begin{eqnarray*}%%%%%%%%%%%%%
    B(X_0,r)\cap \Gamma = \{(x,y)\in B(X_0,r): x\geq x_0 \}, \quad \text{for some } r>0.
\end{eqnarray*}%%----------------------------%%
Fix $h<0$ small enough and let $\phi\in C_c( B(X_0,r))$ be a cut-off function such that $\phi=1$ on $B(X_0,\frac{r}2)$. We choose $(\bb{q},\psi) = (\bb{R}_w\big[D_{h,-\bb{e}_2}^s 	D_{h,\bb{e}_2}^s[w\phi]\big], D_{h,-\bb{e}_2}^s 	D_{h,\bb{e}_2}^s[w\phi])$ in the coupled momentum equation $\eqref{weaksolmom}$ (note that the direction of $\bb{e}_2$ which appears in the quotient is parallel with $B(X_0,r)\cap \partial\Gamma$). Therefore, we can once again close the estimates in the same way as in the periodic domain case to obtain
\begin{eqnarray}\label{ydir}%%%%%%%%%%%%%
    ||D_{h,\bb{e}_2}^s[\Delta w]||_{L^2(0,T; L^2(B(X_0,\frac{r}2)))}\leq C.
\end{eqnarray}%%----------------------------%%
It remains to prove that
\begin{eqnarray*}%%%%%%%%%%%%%
    ||D_{h,\bb{e}_1}^s[w]||_{L^2(0,T; L^2(B(X_0,\frac{r}2)))}\leq C.
\end{eqnarray*}%%----------------------------%%
For this purpose, we can extend $w$ to be defined on $(0,T)\times B(X_0,r)$ and odd with respect to the boundary $x=x_0$
\begin{eqnarray*}%%%%%%%%%%%%%
    \hat{w}(t,x,y):=\begin{cases}
    w(t,x,y),&\quad \text{for } x\geq x_0,\\
    -w(t,x_0+(x_0-x),y),&\quad \text{for } x< x_0.
    \end{cases}
\end{eqnarray*}%%----------------------------%%
This gives us that the clamped boundary condition is satisfied for the following test function (here w.l.o.g. we assume that $\phi$ is also odd with respect to $\{x=x_0\}$)
\begin{eqnarray*}%%%%%%%%%%%%%
    \partial_{\bb{e}_1}D_{h,-\bb{e}_1}^s D_{h,\bb{e}_1}^s[\hat{w}\phi] = D_{h,-\bb{e}_1}^s D_{h,\bb{e}_1}^s[\hat{w}\phi] = 0, \quad \text{on } B(X_0,r)\cap \{x=x_0\},
\end{eqnarray*}%%----------------------------%%
while outside of $B(X_0,r)$ we have that $\hat{w}\phi=0$, so we can test $\eqref{weaksolmom}$ with the pair \\ $(\bb{R}_w\big[D_{h,-\bb{e}_1}^s D_{h,\bb{e}_1}^s[\hat{w}\phi]\big], D_{h,-\bb{e}_1}^s D_{h,\bb{e}_1}^s[\hat{w}\phi])$ and once again obtain 
\begin{eqnarray*}%%%%%%%%%%%%%
    ||D_{h,\bb{e}_2}^s[\Delta w]||_{L^2(0,T; L^2(B(X_0,\frac{r}2)))}\leq C,
\end{eqnarray*}%%----------------------------%%
which combined with $\eqref{ydir}$ finally implies
\begin{eqnarray*}%%%%%%%%%%%%%
    ||\Delta w||_{L^2(0,T; H^s(B(X_0,\frac{r}2)))}\leq C.
\end{eqnarray*}%%----------------------------%%

\noindent
\textbf{Step 3 - boundary estimates (general boundary case)}. Here, we assume that the part of the boundary we are dealing with is a sub-graph of a $C^{0,1}$ function. The idea is to ``straighten out the boundary'' by means of a local coordinate transform. For that purpose, denote by $x',y'$ variables in new coordinate system and $(\bb{e}_1',\bb{e}_2')$ the corresponding unit basis vectors and for $X_0\in \partial\Gamma$ let
\begin{eqnarray*}%%%%%%%%%%%%%
    \Phi:B(0,r_1)\to B(X_0,r_2)
\end{eqnarray*}%%----------------------------%%
be a $C^{0,1}$ local coordinate transform such that 
\begin{eqnarray*}%%%%%%%%%%%%%
\Phi\big(B(0,r_1)\cap\{x'> 0\}\big)= B(X_0,r_2)\cap \Gamma,\qquad 
    \Phi\big(B(0,r_1)\cap \{x=0\}\big)=  B(X_0,r_2)\cap \partial\Gamma.
\end{eqnarray*}%%----------------------------%%
We can express
\begin{eqnarray*}%%%%%%%%%%%%%
    \int_{B(X_0,r_2)\cap \Gamma} \Delta w \Delta \psi ~dx dy= \int_{B(0,r_1)\cap \{x'> 0\}}|\det \nabla \Phi| (\Delta w)\circ \Phi (\Delta \psi)\circ \Phi~ dx'dy' . 
\end{eqnarray*}%%----------------------------%%
Now, we can repeat the calculation from step 2 by choosing $\psi$ to satisfy $(\Delta \psi)\circ \Phi = D_{h,-\bb{e}_2'}^s D_{h,\bb{e}_2'}^s \big[(\Delta w)\circ\Phi \big]$ on $B(0,r_1)$, or equivalently
\begin{eqnarray*}%%%%%%%%%%%%%
    \psi = \Delta^{-1}\Big[\Phi^{-1} \big(D_{h,-\bb{e}_2'}^s D_{h,\bb{e}_2'}^s \big[(\Delta w)\circ\Phi\big]\big)\Big], \quad \text{ on }B(X_0,r_2),
\end{eqnarray*}%%----------------------------%%
where $\Delta^{-1}$ is the inverse Laplace operator with zero Dirichlet boundary data, while the fluid test function is chosen as $\bb{q}=\bb{R}_w[\psi]$. With these test functions, the equation $\eqref{weaksolmom}$ leads to the following estimate
\begin{eqnarray*}%%%%%%%%%%%%%
    &&c\int_{B(0,r_1)\cap \{x'> 0\}} |D_{h,\bb{e}_2'}^s[(\Delta (\phi w))\circ \Phi]|^2 ~ dx'dy'\\
    &&\leq \int_{B(0,r_1)\cap \{x'> 0\}}|\det \nabla \Phi| |D_{h,\bb{e}_2'}^s[(\Delta (\phi w))\circ \Phi]|^2~ dx'dy' \leq C
\end{eqnarray*}%%----------------------------%%
and similarly in the direction $\bb{e}_1'$
\begin{eqnarray*}%%%%%%%%%%%%%
    \int_{B(0,r_1)\cap \{x'> 0\}} |D_{h,\bb{e}_1'}^s[(\Delta (\phi w))\circ \Phi]|^2~ dx'dy' \leq C.
\end{eqnarray*}%%----------------------------%%
Thus, we conclude
\begin{eqnarray*}%%%%%%%%%%%%%
   ||\Delta (\phi w)||_{L^2(0,T; H^s(B(X_0,r_2)\cap \Gamma) )}\leq C||[\Delta (\phi w)] \circ \Phi ||_{L^2(0,T; H^s(B(0,r_1)\cap \{x'> 0\}) )} \leq C
\end{eqnarray*}%%----------------------------%%
because $\Phi$ is a $C^{0,1}$ bijective map. Since $X_0$ was an arbitrary point and boundary $\partial\Gamma$ can be covered by a finite number of balls $B(X_i,\frac{r_1}2)$ with non-empty intersection such that $X_i\in \partial\Gamma$ and $B(X_i,r_1)\cap \partial\Gamma$ is a sub-graph of a $C^{0,1}$ function, the proof is finished.

\subsection{Part 2: relative energy inequality - Theorem $\ref{main2}(2)$}
We start by choosing  $(\bb{q},\psi)=(\bb{u}_2, \partial_t w_2)$ in the coupled momentum equation $\eqref{weaksolmom}$ on\footnote{Note that the coupled momentum equation $\eqref{weaksolmom}$ and the continuity equation $\eqref{cont}$ hold on $(0,t)$, for any $t\in (0,T]$, due to the weak continuity in time of $\rho_1$ and $\rho_1 \bb{u}_1$ and $\partial_t w_1$.} $(0,t)$ to obtain
\begin{eqnarray}%%%%%%%%%%%%%
&&\int_{\Omega^{w_1}(t)} (\rho_1 \bb{u}_1 \cdot \bb{u}_2)(t)+ \int_\Gamma (\partial_t w_1 \partial_t w_2)(t) \nonumber \\
&&= \int_{\Omega^{w_1}(0)} (\rho_1 \bb{u}_1 \cdot \bb{u}_2)(0) + \int_\Gamma (\partial_t w_1 \partial_t w_2)(0)  + \int_{Q_t^{w_1}}\rho_1 \bb{u}_1\cdot\partial_t \bb{u}_2\nonumber \\
&&\quad+ \int_{Q_t^{w_1}} ( \rho_1\bb{u}_1 \otimes \bb{u}_1) :\nabla \bb{u}_2  + \int_{Q_t^{w_1}}  \rho_1^\gamma (\nabla \cdot \bb{u}_2) 
- \underbrace{\int_{Q_t^{w_1}}  \mathbb{S}(\nabla \bb{u}_1):\nabla\bb{u}_2}_{=\int_{Q_t^{w_1}}  \mathbb{S}(\nabla\bb{u}_2):\nabla \bb{u}_1} \nonumber\\
&&\quad + \int_{\Gamma_t} \partial_t w_1 \partial_t^2 w_2- \int_{\Gamma_t} \Delta w_1 \Delta \partial_t w_2-
\alpha \bint_{\Gamma_t}\partial_t \nabla w_1 \cdot \partial_t \nabla w_2. \label{ws1}
\end{eqnarray}%%----------------------------%%
Next, in the continuity equation $\eqref{cont}$ on $(0,t)$, choose $\varphi = \frac{1}{2}|\bb{u}_2|^2$ to obtain
\begin{eqnarray}%%%%%%%%%%%%%
&&\frac{1}{2} \int_{\Omega^{w_1}(t)} (\rho_1|\bb{u}_2|^2)(t)=\frac{1}{2} \int_{\Omega^{w_1}(0)} (\rho_1|\bb{u}_2|^2)(0)+ \int_{Q_t^{w_1}} \rho_1\partial_t \bb{u}_2 \cdot \bb{u}_2 + \int_{Q_t^{w_1}} \rho_1( \bb{u}_1\cdot \nabla \bb{u}_2)\cdot \bb{u}_2 , \quad\quad\label{ws2}
\end{eqnarray}%%----------------------------%%
and then choose $\varphi = \frac{\gamma}{\gamma-1} \rho_2^{\gamma-1}$
\begin{align}\label{ws3}%%%%%%%%%%%%%
&\frac{\gamma}{\gamma-1}\int_{\Omega^{w_1}(t)} (\rho_1 \rho_2^{\gamma-1})(t)\nonumber\\
&= \frac{\gamma}{\gamma-1}\int_{\Omega^{w_1}(0)} (\rho_1 \rho_2^{\gamma-1})(0)+\frac{\gamma}{\gamma-1}\int_{Q_t^{w_1}}\rho_1\partial_t (\rho_2^{\gamma-1}) + \frac{\gamma}{\gamma-1}\int_{Q_t^{w_1}}  \rho_1\bb{u}_1\cdot\nabla (\rho_2^{\gamma-1}) .
\end{align}%%----------------------------%%
Now, by Theorem $\ref{main2}(1)$, one has that $w\in L^2(0,T;C^{0,1}(\Gamma))$ which means that $\Omega^w(t)$ is Lipschitz for a.a. $t\in (0,T)$. This allows us to use the divergence theorem to obtain
\begin{eqnarray*}%%%%%%%%%%%%%
	\frac{\gamma}{\gamma-1}\int_{Q_t^{w_1}} \rho_2 \bb{u}_2 \cdot\nabla (\rho_2^{\gamma-1}) = \int_{Q_t^{w_1}}  \bb{u}_2 \cdot\nabla (\rho_2^\gamma) = - \int_{Q_t^{w_1}} \rho_2^\gamma (\nabla \cdot \bb{u}_2) + \int_{\Gamma_t^{w_1}}  \rho_2^\gamma \partial_t w_2\nu^{w_1} \cdot\bb{e}_3 ,
\end{eqnarray*}%%----------------------------%%
and the Raynolds transport theorem to obtain
\begin{eqnarray*}%%%%%%%%%%%%%
	\int_0^t\frac{d}{dt} \int_{\Omega^{w_1}(t)} \rho_2^\gamma = \frac{\gamma}{\gamma-1}\int_{Q_t^{w_1}}\rho_2\partial_t (\rho_2^{\gamma-1})+ \int_{\Gamma_t^{w_1}} \rho_2^{\gamma}\partial_t w_1 \nu^{w^1} \cdot \bb{e}_3,
\end{eqnarray*}%%----------------------------%%
which together imply
\begin{align}%%%%%%%%%%%%%
& \int_{\Omega^{w_1}(t)} \rho_2^{\gamma}(t)= \int_{\Omega^{w_1}(0)}  \rho_2^{\gamma}(0)+\int_{\Gamma_t^{w_1}} \rho_2^\gamma(\partial_t w_1 - \partial_t w_2) \nu^{w_1} \cdot\bb{e}_3 \nonumber \\
& + \frac{\gamma}{\gamma-1}\int_{Q_t^{w_1}} \Big[\rho_2 \bb{U}_2 \cdot \nabla (\rho_2^{\gamma-1}) + \rho_2 \partial_t (\rho_2^{\gamma-1})\Big] + \int_{Q_t^{w_1}} \rho_2^{\gamma} (\nabla \cdot \bb{u}_2). \label{ws5}
\end{align}%%----------------------------%%
Finally, denoting $a\pm b=a+b-b$, the last three terms on the RHS of $\eqref{ws1}$ can be transformed as
\begin{eqnarray*}%%%%%%%%%%%%%
	&&\int_{\Gamma_t} \Delta w_1 \partial_t \Delta w_2\pm \frac{1}{2}\Big( \int_\Gamma |\Delta w_2|^2\Big)\Big|_0^t= \int_{\Gamma_t} w_1 \partial_t \Delta^2 w_2 - \frac12\int_0^t \frac{d}{dt}\int_\Gamma |\Delta w_2|^2  + \frac{1}{2}\Big( \int_\Gamma |\Delta w_2|^2\Big)\Big|_0^t \\
	&&= - \int_{\Gamma_t} \partial_t w_1 \Delta^2 w_2 + \Big(\int_\Gamma w_1 \Delta^2 w_2\Big)\Big|_0^t - \int_{\Gamma_t} \partial_t\Delta w_2 \Delta w_2+ \frac12\Big( \int_\Gamma |\Delta w_2|^2\Big)\Big|_0^t  \\
	&&= \int_{\Gamma_t} (\partial_t w_2 - \partial_t w_1) \Delta^2 w_2 +  \Big(\int_\Gamma \Delta w_1 \Delta w_2\Big)\Big|_0^t+ \frac12\Big( \int_\Gamma |\Delta w_2|^2\Big)\Big|_0^t,
\end{eqnarray*}%%----------------------------%%	
and
\begin{eqnarray*}%%%%%%%%%%%%%	
	&&- \int_{\Gamma_t} \partial_t w \partial_t^2 w_2 \pm \Big(\frac{1}{2} \int_\Gamma |\partial_t w_2|^2\Big)\Big|_0^t = \int_{\Gamma_t} (\partial_t w_2 - \partial_t w_1) \partial_t^2 w_2 + \Big(\frac{1}{2} \int_\Gamma |\partial_t w_2|^2\Big)\Big|_0^t~ ,\\
	&&-\bint_{\Gamma_t}\partial_t \nabla w_1 \cdot \partial_t\nabla w_2 \pm \int_{\Gamma_t} |\partial_t \nabla w_1|^2 =	\int_{\Gamma} \partial_t \Delta w_2 (\partial_t  w_1 - \partial_t w_2)+\int_{\Gamma_t} |\partial_t \nabla w_1|^2~ ,
\end{eqnarray*}%%----------------------------%%
so we sum $\eqref{energy}-\eqref{ws1}+\eqref{ws2}-\eqref{ws3}-\eqref{ws5}$, which then gives us $\eqref{relen}$.

\section{Proof of Theorem $\ref{main1}$}\label{proofofmain2}
When working with free boundary problems, it is common to transform the coordinate system in order to define the problem on a fixed domain. Also, in this proof, we will be working with two solutions which are in general defined on different domains, so transforming them onto a common domain will allow us to compare them. 

\subsection{Part 1: change of coordinates}\label{sec3}
We start by introducing the fixed domain
\begin{eqnarray*}%%%%%%%%%%%%%
	\Omega:= \Gamma \times (-1,0),
\end{eqnarray*}%%----------------------------%%
and the family of so called arbitrary Lagrangian-Eulerian (ALE) mappings 
\begin{eqnarray}%%%%%%%%%%%%%
A_w(t):&&\Omega\to\Omega^w(t)\nonumber\\
&&(X,z)\mapsto(X, (z+1)w(t,X)+z). \label{ALE}
\end{eqnarray}%%----------------------------%%
For an arbitrary (vector or scalar) function $\bb{f}$ defined on $\Omega^w (t)$, corresponding to $A_w$, we denote the \textbf{pull-back}, \textbf{transformed gradient} and \textbf{transformed divergence} of $\bb{f}$ on $\Omega$, respectively, as
\begin{eqnarray*}%%%%%%%%%%%%%
	&&\bb{f}^w:=\bb{f} \circ A_w,\\
	&&\nabla^w \mathbf{f}^w:=(\nabla \mathbf{f})^w = \nabla \mathbf{f}^w (\nabla A_w)^{-1}\circ A_w, \\
	&&\nabla^w \cdot \bb{f}^w := \text{Tr}(\nabla^w \mathbf{f}^w).
\end{eqnarray*}%%----------------------------%%
where
\begin{eqnarray*}
(\nabla A_w^{-1})\circ A_w =\begin{bmatrix}
    1  & 0 & 0\\
    0  & 1 & 0\\
   -\ddfrac{z+1}{w+1}\partial_x w & -\ddfrac{z+1}{w+1}\partial_y w & \ddfrac{1}{w+1}
\end{bmatrix}.
\end{eqnarray*}
The \textbf{Jacobian} and the \textbf{ALE velocity} of the mapping $A_w$, respectively, read
\begin{eqnarray*}%%%%%%%%%%%%%
	J&:=&\text{det}\nabla A_w = 1+w,\\
	\bb{w}&:=& \frac{d}{dt} A_w = (z+1)\partial_t w \bb{e}_3,
\end{eqnarray*}%%----------------------------%%
and the time-space cylinder is denoted as
\begin{eqnarray*}%%%%%%%%%%%%%
	Q_T:= (0,T)\times\Omega.
\end{eqnarray*}%%----------------------------%%
Since we will work with pairs of solutions, the following notation will be useful. For $i=1,2$ and a given (scalar or vector function) $\bb{f}$, denote
\begin{eqnarray*}%%%%%%%%%%%%%
\nabla_i \bb{f}:= \nabla^{w_i} \bb{f},\quad \nabla_i \cdot \bb{f}:= \nabla^{w_i}\cdot \bb{f},\quad A_i:=A_{w_i},\quad S^{i}:= S^{w^i},\quad \nu^i := \nu^{w^i},
\end{eqnarray*}%%----------------------------%%
and the pull-backs
\begin{eqnarray*}%%%%%%%%%%%%%
	r_i&:=&\rho_i^{w_i} = \rho_i \circ A_{i}, \\
	\bb{U}_i&:=&\bb{u}_i^{w_i} = \bb{u}_i \circ A_{i}.
\end{eqnarray*}%%----------------------------%%
Finally, the subscript for the Jacobian will be omitted
\begin{eqnarray*}%%%%%%%%%%%%%
	J:=J_1 = 1+w_1,
\end{eqnarray*}%%----------------------------%%
since it will appear frequently.

\subsubsection{The relative entropy and relative energy inequality on the fixed reference domain}
Let $(\rho_1,\bb{u}_1,w_1)$ be a finite energy weak solution in the sense of Definition $\ref{weaksolution}$ and let $r_2,\bb{U}_2 \in C^\infty([0,T]\times \overline{\Omega})$ and $w_2 \in C_0^\infty(\Gamma_T)$ such that ${\bb{U}_2}_{|\Sigma}=0$ and ${\bb{U}_2}_{|\Gamma\times\{0\}} = \partial_t w_2\bb{e}_3$. For\footnote{Here, the weak solution on the fixed reference domain is not defined explicitly since it will not be required in the proof, but this definition can be found in \cite[section 2.3]{trwa3}.} $(r_1,\bb{U}_1,w_1) = (\rho_1 \circ A_1,\bb{u}_1 \circ A_1, w_1)$, the \textbf{relative entropy on the fixed reference domain} with respect to $(r_2,\bb{U}_2,w_2)$ will be also denoted by $\mathcal{E}$ and it reads as \\
\begin{eqnarray}%%%%%%%%%%%%%
\mathcal{E}\Big((r_1,\bb{U}_1,w_1)\Big|(r_2,\bb{U}_2,w_2) \Big)(t) &=& \frac{1}{2} \int_\Omega (J r_1 |\bb{U}_1 - \bb{U}_2|^2)(t) + \frac{1}{\gamma-1}\int_{\Omega}J \big( r_1^\gamma - \gamma r_2^{\gamma-1}(r_1 - r_2) - r_2^{\gamma} \big)(t) \nonumber \\
&&\quad +\frac{1}{2}\int_\Gamma |\partial_t w_1 - \partial_t w_2|^2(t) + \frac{1}{2}\int_\Gamma |\Delta w_1 - \Delta w_2|^2(t). \nonumber \\\nonumber
\end{eqnarray}%%----------------------------%%
By using Theorem $\ref{main2}(2)$ and the identities
\begin{eqnarray}%%%%%%%%%%%%%
&&\int_{\Gamma^{w_1}(t)} (\partial_t w_2 - \partial_t w_1)\rho_2^\gamma I\nu^1\cdot \bb{e}_3 =\int_{\Gamma} (\partial_t w_2 - \partial_t w_1)r_2^\gamma \underbrace{S_1 I \nu^1\cdot \bb{e}_3}_{=1}, \nonumber \\
&&\ddfrac{d}{dt} (\bb{q} \circ A_i)= (\partial_t \bb{q})\circ A_i + \bb{w}_i \cdot \nabla_i (\bb{q} \circ A_i),\quad \text{ on}\quad \Omega,  \label{aleid} 
\end{eqnarray}%%-------
where $\bb{q}$ is a sufficiently regular function defined on $Q_T^{w_i}$, one can easily see that the triplet $(r_1,\bb{U}_1,w_1)$ satisfies the \textbf{relative energy inequality on the fixed reference domain} with respect to $(r_2,\bb{U}_2,w_2)$ of the form \\
\begin{eqnarray}%%%%%%%%%%%%%
&&\mathcal{E}\Big((r_1,\bb{U}_1,w_1)\Big|(r_2,\bb{U}_2,w_2) \Big)(t) + \int_{Q_t} J \mathbb{S} (\nabla_1 \bb{U}_1 - \nabla_1 \bb{U}_2):(\nabla_1\bb{U}_1 - \nabla_1 \bb{U}_2)+ \alpha\int_0^t\int_{\Gamma}  |\partial_t \nabla w_1 - \partial_t \nabla w_2|^2 \nonumber\\
&&\leq \mathcal{E}\Big((r_1,\bb{U}_1,w_1)\Big|(r_2,\bb{U}_2,w_2) \Big)(0) +\int_0^t\mathcal{R}\Big( r_1,\bb{U}_1,w, r_2,\bb{U}_2,w_2\Big)\label{relenfixed} \\ \nonumber
\end{eqnarray}%%----------------------------%%
where the \textbf{remainder term} reads \\
\begin{align}%%%%%%%%%%%%%
&\mathcal{R}\Big( r_1,\bb{U}_1,w_1, r_2,\bb{U}_2,w_2\Big) := \int_{\Omega} J \mathbb{S}( \nabla_1 \bb{U}_2):(\nabla_1\bb{U}_2 - \nabla_1 \bb{U}_1)+\int_{\Omega}J r_1 (\partial_t \bb{U}_2 + \bb{U}_1\cdot \nabla_1 \bb{U}_2)\cdot(\bb{U}_2 - \bb{U}_1) \nonumber\\
&\quad- \int_{\Omega} J r_1 \big[\bb{w}_1 \cdot \nabla_1 \bb{U}_2\big] \cdot (\bb{U}_2-\bb{U}_1)+ \frac{\gamma}{\gamma-1} \Big[\int_{\Omega} J (r_2\bb{U}_2 - r_1 \bb{U}_1)\cdot \nabla_1 (r_2^{\gamma-1})  \nonumber\\
&\quad+ J(r_1 - r_2)\big( \partial_t(r_2^{\gamma-1}) -\bb{w}_1 \cdot \nabla_1(r_2^{\gamma-1})\big)\Big]
+ \int_{\Omega} J(r_2^\gamma - r_1^\gamma)(\nabla_1 \cdot \bb{U}_2) \nonumber \\
&\quad+ \int_{\Gamma}(\partial_t w_2 - \partial_t w_1) r_2^\gamma -\int_{\Gamma} (\partial_t w_2 - \partial_t w_1) \partial_t^2 w_2 -\int_{\Gamma} (\partial_t w_2 - \partial_t w_1) \Delta^2 w_2+\alpha\int_{\Gamma} \partial_t \Delta w_2 (\partial_t  w_1 - \partial_t w_2). \label{caprfixed} 
\end{align}%%----------------------------%%
\subsection{Part 2: carrying out the proof on the fixed reference domain}
\begin{figure}[H]
	\centering\includegraphics[scale=0.48]{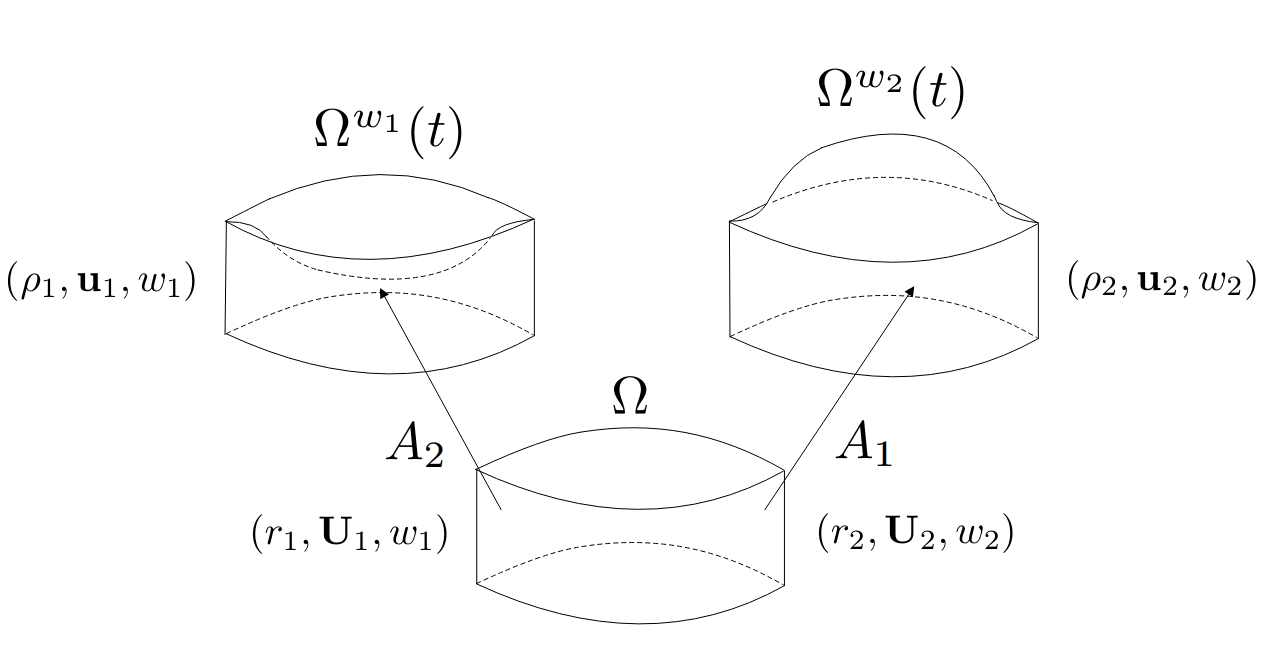}
	\caption{The weak and the strong solution on both physical and fixed reference domain coordinates.}
	\label{domains}
\end{figure}

Here, we compare (see figure $\ref{domains}$) a finite energy weak solution $(\rho_1,\bb{u}_1,w_1)$ and a strong solution $(\rho_2,\bb{u}_2,w_2)$ belonging to the regularity class specified in Theorem $\ref{main1}$, which emanate from the same initial data. 

The goal is to prove that 
\begin{eqnarray*}%%%%%%%%%%%%%
\mathcal{E}\Big((r_1,\bb{U}_1,w_1)\Big|(r_2,\bb{U}_2,w_2) \Big)(t)\leq C\int_0^t h(\tau) 
\mathcal{E}\Big((r_1,\bb{U}_1,w_1)\Big|(r_2,\bb{U}_2,w_2) \Big)(\tau) d\tau, \quad h \in L^1(0,T),
\end{eqnarray*}%%----------------------------%%
for a.a. $t\in(0,T_2)$, then the desired result will follow immediately by Gronwall's lemma.\\

${}$\\
\noindent
\textbf{\textit{Step 1: transforming the remainder term} $\mathcal{R}\Big( r_1,\bb{U}_1,w_1, r_2,\bb{U}_2,w_2\Big) $}\\

First, one can see that $(r_2,\bb{U}_2,w_2)$ satisfies the following transformed problem in the strong sense
\begin{align}\label{strongfixed}%%%%%%%%%%%%%
\begin{cases}
\partial_t^2 w_2 + \Delta^2 w_2-\alpha\partial_t \Delta w_2 = -S^{2}\big[(\mathbb{S}(\nabla_2 \bb{U}_2) - r_2^\gamma I ) \nu^{2}\big] \cdot \bb{e}_3,\quad \text{ on  }~ \Gamma_{T_2}, \\[2.5mm]
\partial_t r_2 - \bb{w}_2\cdot \nabla_2 r_2 + \nabla_2 \cdot (r_2\bb{U}_2 ) = 0,\quad \text{ in  }~ Q_{T_2}, \\[2.5mm]
\partial_t \bb{U}_2 - \bb{w}_2\cdot \nabla_2\bb{U}_2 + \bb{U}_2
\cdot \nabla_2 \bb{U}_2 = - \frac{1}{r_2}\nabla_2 (r_2^\gamma) + \frac{1}{r_2} \nabla_2\cdot \mathbb{S}(\nabla_2\bb{U}_2), \quad \text{ in  }~Q_{T_2},
\end{cases}
\end{align}%%----------------------------%%
by using $\eqref{aleid}$, and that
\begin{eqnarray*}%%%%%%%%%%%%%
	&&r_2 \in L^2(0,T_2; W^{1,q}(\Omega)) \cap H^1(0,T_2; L^q(\Omega)), \\
	&&\bb{U}_2 \in L^2(0,T_2; W^{2,q}(\Omega)) \cap H^1(0,T_2; L^q(\Omega)),\\
	&&w_2 \in L^2(0,T_2; H^4(\Gamma)) \cap H^2(0,T_2; L^2(\Gamma)),\\
	&& 0< \inf\limits_{Q_{T_2}} r_2\leq \sup\limits_{Q_{T_2}} r_2 <\infty,
\end{eqnarray*}%%----------------------------%%
with $q = q(d,\alpha,\gamma)$ being specified in Theorem $\ref{main1}$. Using this regularity, we test $\eqref{relenfixed}$ with $(r_2,\bb{U}_2,w_2)$ (by the density argument), and then transform the remainder term given in $\eqref{caprfixed}$ by  $\eqref{strongfixed}$, as follows. \\

\noindent
$\quad\quad$\textbf{The 2nd term on RHS of }$\eqref{caprfixed}$:
\begin{eqnarray}%%%%%%%%%%%%%
&&\int_{Q_t} Jr_1 (\partial_t \bb{U}_2 + \bb{U}_1 \cdot \nabla_1 \bb{U}_2)(\bb{U}_2 - \bb{U}_1)\nonumber \\
&&= \int_{Q_t} Jr_1 (\partial_t \bb{U}_2 + \bb{U}_2 \cdot \nabla_2 \bb{U}_2)(\bb{U}_2 - \bb{U}_1) + \int_{Q_t} Jr_1 \big[\bb{U}_2\cdot (\nabla_1 - \nabla_2)\bb{U}_2 \big]\cdot(\bb{U}_2 - \bb{U}_1)\nonumber \\
&&\quad+ \int_{Q_t} Jr_1 \big[(\bb{U}_1 - \bb{U}_2) \cdot \nabla_1 \bb{U}_2 \big]\cdot(\bb{U}_2 - \bb{U}_1) \pm \int_{Q_t}Jr_1(\bb{w}_2\cdot \nabla_2\bb{U}_2)(\bb{U}_2 - \bb{U}_1) \nonumber\\
&&\stackrel{\eqref{strongfixed}_3}{=} -\int_{Q_t} J\frac{r_1}{r_2}\nabla_2(r_2^\gamma)(\bb{U}_2-\bb{U}_1) +\int_{Q_t} J \frac{1}{r_2} \big[\nabla_2\cdot \mathbb{S}(\nabla_2\bb{U}_2)\big]\cdot (\bb{U}_2 - \bb{U}_1) \nonumber\\
&&\quad+ \int_{Q_t} Jr_1 \big[\bb{U}_2\cdot (\nabla_1 - \nabla_2)\bb{U}_2\big]\cdot(\bb{U}_2 - \bb{U}_1) + \int_{Q_t} Jr_1 \big[(\bb{U}_1 - \bb{U}_2) \cdot \nabla_1 \bb{U}_2 \big]\cdot(\bb{U}_2 - \bb{U}_1)\nonumber\\
&&\quad+ \int_{Q_t}Jr_1(\bb{w}_2\cdot \nabla_2\bb{U}_2)(\bb{U}_2 - \bb{U}_1). \label{bound1}
\end{eqnarray}%%----------------------------%%

\textbf{The 2nd term on RHS of }$\eqref{bound1}$:
\begin{eqnarray}%%%%%%%%%%%%%
&&\int_{Q_t} J \frac{1}{r_2} r_1\big[\nabla_2\cdot \mathbb{S}(\nabla_2\bb{U}_2)\big]\cdot(\bb{U}_2 - \bb{U}_1)\nonumber \\
&&=\int_{Q_t}J\frac{1}{r_2}(r_1 - r_2)\big[\nabla_2 \cdot \mathbb{S}(\nabla_2\bb{U}_2)\big]\cdot (\bb{U}_2 - \bb{U}_1)+ \int_{Q_t}(J-J_2)\big[\nabla_2 \cdot \mathbb{S}(\nabla_2\bb{U}_2)\big]\cdot (\bb{U}_2 - \bb{U}_1) \nonumber\\
&&\quad+\int_{Q_t}J_2\big[\nabla_2 \cdot \mathbb{S}(\nabla_2\bb{U}_2)\big]\cdot (\bb{U}_2 - \bb{U}_1), \label{bound2}
\end{eqnarray}%%----------------------------%%
where the last term can be transformed by partial integration on the physical domain
\begin{eqnarray}%%%%%%%%%%%%%
	&&\int_{Q_t}J_2\nabla_2 \cdot \mathbb{S}(\nabla_2\bb{U}_2)(\bb{U}_2 - \bb{U}_1) = \int_{0}^t \int_{\Omega^{w_2}(t)} \nabla \cdot \mathbb{S}(\nabla\bb{u}_2)(\bb{u}_2 - \bb{U}_1 \circ A_2^{-1})\nonumber \\
	&&= -\int_{0}^t \int_{\Omega^{w_2}(t)} \mathbb{S}(\nabla\bb{U}_2):(\nabla\bb{u}_2 -\nabla \bb{U}_1\circ A_2^{-1}) + \int_0^t\int_{\Gamma^{w_2}(t)} \big[(\mathbb{S}(\nabla \bb{u}_2))\nu^2\big] \cdot \bb{e}_3(\partial_t w_2 - \partial_t w_1) \nonumber  \\	
	 &&=-\int_{Q_t}J_2 \mathbb{S}(\nabla_2\bb{U}_2):(\nabla_2\bb{U}_2 - \nabla_2\bb{U}_1) + \int_{\Gamma_t}S_2\big[(\mathbb{S}(\nabla_2 \bb{U}_2))\nu^2\big] \cdot \bb{e}_3(\partial_t w_2 - \partial_t w_1). \label{partial}
\end{eqnarray}%%----------------------------%%

${}$

\textbf{The 4th term on RHS of} $\eqref{caprfixed}$ and \textbf{the 1st term on RHS of} $\eqref{bound1}$:
\begin{eqnarray}%%%%%%%%%%%%%
&&\frac{\gamma}{\gamma-1} \Big[\int_{Q_t} J (r_2\bb{U}_2 - r_1 \bb{U}_1) \cdot \nabla_1 (r_2^{\gamma-1}) + J(r_1 - r_2)\big[ \partial_t(r_2^{\gamma-1}) -\bb{w}_1 \cdot \nabla_1(r_2^{\gamma-1})\big] \Big] \nonumber\\
&&- \int_{Q_t}\underbrace{ J \frac{r_1}{r_2}\nabla_2(r_2^\gamma)}_{=\frac{\gamma}{\gamma-1}J r_1 \nabla_2 (r_2^{\gamma-1}) }(\bb{U}_2-\bb{U}_1)\nonumber \\
&& = \frac{\gamma}{\gamma-1} \int_{Q_t} J (r_2 - r_1 )\Big[ \bb{U}_2\cdot\nabla_1 (r_2^{\gamma-1})+\partial_t (r_2^{\gamma-1}) -\bb{w}_1\cdot \nabla_1(r_2^{\gamma-1}) \Big] \nonumber \\
&&\quad+ \int_{Q_t} r_1(\nabla_1 - \nabla_2)(r_2^{\gamma-1})(\bb{U}_2 - \bb{U}_1)\pm \frac{\gamma}{\gamma-1} \int_{Q_t} J (r_2 - r_1 )\Big[ \bb{U}_2\cdot \nabla_2 (r_2^{\gamma-1}) -\bb{w}_2 \cdot \nabla_2(r_2^{\gamma-1}) \Big] \nonumber\\
&&= \frac{\gamma}{\gamma-1} \int_{Q_t} J (r_2 - r_1 )\Big[\underbrace{\partial_t (r_2^{\gamma-1})+ \bb{U}_2\cdot \nabla_2 (r_2^{\gamma-1}) -\bb{w}_2 \cdot \nabla_1(r_2^{\gamma-1})}_{=-(\gamma-1)(\nabla_2 \cdot \bb{U}_2) r_2^{\gamma-1}, ~~\text{by } \eqref{strongfixed}_2 }\Big] \nonumber \\
&&\quad +\frac{\gamma}{\gamma-1} \int_{Q_t} J(r_2 - r_1 )\Big[ \bb{U}_2 \cdot (\nabla_1 - \nabla_2)(r_2^{\gamma-1}) - (\bb{w}_1 - \bb{w}_2)\cdot\nabla_1(r_2^{\gamma-1}) - \bb{w}_2\cdot (\nabla_1 - \nabla_2)(r_2^{\gamma-1}) \Big]\nonumber \\                             
&&\quad+ \int_{Q_t}J r_1 (\nabla_1 - \nabla_2)(r_2^{\gamma-1})\cdot (\bb{U}_2 - \bb{U}_1). \label{bound3}
\end{eqnarray}%%----------------------------%%

\textbf{The 5th term on RHS of }$\eqref{caprfixed}$ and \textbf{the 1st term on RHS of }$\eqref{bound3}$:
\begin{eqnarray}%%%%%%%%%%%%%
&&\int_{Q_t} J(r_2^\gamma - r_1^\gamma)(\nabla_1 \cdot \bb{U}_2) - \gamma \int_{Q_t}J(r_2 - r_1)(\nabla_2 \cdot \bb{U}_2) r_2^{\gamma-1} \nonumber\\
&&=\int_{Q_t} J\Big[ r_2^\gamma - \gamma r_2^{\gamma-1}(r_2 - r_1) - r_1^\gamma\Big](\nabla_1 \cdot \bb{U}_2)+ \gamma\int_{Q_t} J (r_2 - r_1)r_2^{\gamma-1}(\nabla_1 - \nabla_2) \cdot \bb{U}_2.\label{bound5}
\end{eqnarray}%%----------------------------%%

\textbf{The plate equation }$\eqref{strongfixed}_1 \times (\partial_t w_2 - \partial_t w_1)$:
\begin{eqnarray}%%%%%%%%%%%%
&&\int_{\Gamma_t} (\partial_t w_2 - \partial_t w_1) \partial_t^2 w_2 +\int_{\Gamma_t} (\partial_t w_2 - \partial_t w_1) \Delta^2 w_2 - \alpha \int_{\Gamma_t} \partial_t \Delta w_2 (\partial_t w_2 - \partial_t w_1) \nonumber\\
&& =\int_{\Gamma_t} (\partial_t w_2 - \partial_t w_1) r_2^\gamma \underbrace{S^2 I \nu^2 \cdot \bb{e}_3}_{=1} - \int_{\Gamma_t}(\partial_t w_2 - \partial_t w_1)S_2\big[(\mathbb{S}(\nabla_2 \bb{U}_2))\nu^2\big] \cdot \bb{e}_3.\quad \label{bound4}
\end{eqnarray}%%----------------------------%%

 Now, the remainder term given in $\eqref{caprfixed}$ is transformed by using the identities $\eqref{bound1},\eqref{bound2},\eqref{partial},\eqref{bound3}$,$\eqref{bound5}$ and $\eqref{bound4}$ (in that order), which then by $\eqref{relenfixed}$ gives us:
\begin{eqnarray}%%%%%%%%%%%%%
&&\mathcal{E}\Big((r_1,\bb{U}_1,w_1)\Big|(r_2,\bb{U}_2,w_2) \Big)(t) + \int_{Q_t} J |\mathbb{S} (\nabla_1 \bb{U}_1 - \nabla_1 \bb{U}_2)|^2 \nonumber\\[2mm]
&&\leq \int_{Q_t} J \mathbb{S}( \nabla_1 \bb{U}_2):(\nabla_1\bb{U}_2 - \nabla_1 \bb{U}_1)  -\int_{Q_t}J_2 \mathbb{S}(\nabla_2\bb{U}_2):(\nabla_2\bb{U}_2 - \nabla_2\bb{U}_1) \nonumber\\[2mm]
&&\quad\int_{Q_t}J\frac{1}{r_2}(r_1 - r_2)\big[\nabla_2 \cdot \mathbb{S}(\nabla_2\bb{U}_2)\big]\cdot(\bb{U}_2 - \bb{U}_1)+ \int_{Q_t}(J-J_2)\big[\nabla_2 \cdot \mathbb{S}(\nabla_2\bb{U}_2)\big]\cdot(\bb{U}_2 - \bb{U}_1)\nonumber\\[2mm]
&&\quad+ \frac{\gamma}{\gamma-1} \int_{Q_t} J(r_2 - r_1 )\Big[ \bb{U}_2 \cdot (\nabla_1 - \nabla_2)(r_2^{\gamma-1}) - (\bb{w}_1 - \bb{w}_2)\cdot\nabla_1(r_2^{\gamma-1}) - \bb{w}_2\cdot(\nabla_1 - \nabla_2)(r_2^{\gamma-1}) \Big]\nonumber \\[2mm]
&&\quad+ \int_{Q_t} Jr_1 \big[ \bb{U}_2\cdot (\nabla_1 - \nabla_2)\bb{U}_2\big]\cdot(\bb{U}_2 - \bb{U}_1) + \int_{Q_t} Jr_1 \big[(\bb{U}_1 - \bb{U}_2) \cdot \nabla_1 \bb{U}_2\big]\cdot(\bb{U}_2 - \bb{U}_1)\nonumber\\[2mm]
&&\quad+ \int_{Q_t}J r_1(\nabla_1 - \nabla_2)(r_2^{\gamma-1})\cdot(\bb{U}_2 - \bb{U}_1)+ \gamma \int_{Q_t} J (r_2 - r_1) r_2^{\gamma-1}(\nabla_1 - \nabla_2) \cdot \bb{U}_2 \nonumber \\[2mm]
&&\quad+ \int_{Q_t} J\big[ r_2^\gamma - \gamma r_2^{\gamma-1}(r_2 - r_1) - r_1^\gamma\big](\nabla_1 \cdot \bb{U}_2)+\int_{Q_t} Jr_1\Big[\bb{w}_1 \cdot \nabla_1 \bb{U}_2 -\bb{w}_2\cdot \nabla_2 \bb{U}_2\Big]\cdot (\bb{U}_1 - \bb{U}_2). \nonumber \\ \label{estimate}
\end{eqnarray}%%----------------------------%%
The aim is to bound the right-hand side by 
\begin{eqnarray}\label{provethis}%%%%%%%%%%%%%
	C(\delta)\int_0^t h 
	\mathcal{E}\Big((r_1,\bb{U}_1,w_1)\Big|(r_2,\bb{U}_2,w_2) \Big)+ \delta \int_{Q_t} J |\mathbb{S} (\nabla_1 \bb{U}_1 - \nabla_1 \bb{U}_2)|^2 ,
\end{eqnarray}%%----------------------------%%
for some $h\in L^1(0,T)$ and a small $\delta$. This will be done in step 3.\\ 

\noindent
\textbf{\textit{Step 2: studying the difference $(r_1-r_2)$}}\\ 

Denote
\begin{eqnarray*}%%%%%%%%%%%%%
f(x,y):= x^{\gamma} - \gamma y^{\gamma -1 }(x -y) - y^\gamma. \quad x,y\geq 0.	
\end{eqnarray*}%%----------------------------%%
It is easy to see that for any $C>0$, $f(x,C)$ is convex in $x$ and $f(C,C) = \partial_x f(C,C) = 0$. This means that for a fixed $C>0$, one can choose $c = c(C)>0$ so that
\begin{eqnarray*}%%%%%%%%%%%%%
	f(x,C) \geq c|x-C|^{2\land \gamma}, \quad  x\in \Big[\frac{C}{2}, 2C\Big].
\end{eqnarray*}%%----------------------------%%
where $2\land \gamma = \min\{2, \gamma\}$. Since $0<c_{r_2}=\inf\limits_{Q_{T_2}} r_2 \leq \sup\limits_{Q_{T_2}} r_2 = C_{r_2}$, there is a constant $c=c(c_{r_2},C_{r_2},\gamma)$ for which
\begin{eqnarray}%%%%%%%%%%%%%
f(x,r_2) &\geq& c|x - r_2|^{2\land \gamma}, \quad \text{for } \frac{r_2}{2} \leq x \leq 2r_2,\label{f1}\\
f(x,r_2) &\geq& c(1+x^{2\land \gamma}), \quad \text{for }x \in \mathbb{R}_0^+ \setminus \Big[\frac{r_2}{2}, 2r_2\Big].\label{f2}
\end{eqnarray}%%----------------------------

\noindent
Now, for $0 \leq r_1(t) \leq \frac{r_2(t)}{2}$, one has
\begin{eqnarray*}%%%%%%%%%%%%%
&&||r_1(t) - r_2(t)||_{L^{2 \land \gamma}\big(\big\{ 0 \leq r_1(t) \leq \frac{r_2(t)}{2} \big\}\big)}^2 \leq \int_{\big\{ 0 \leq r_1(t) \leq \frac{r_2(t)}{2} \big\}} c \\
&&\leq \int_{\big\{ 0 \leq r_1(t) \leq \frac{r_2(t)}{2} \big\}}c(1+r_1^\gamma(t)) \leq C\mathcal{E}\Big((r_1,\bb{U}_1,w_1)|(r_2,\bb{U}_2,w_2) \Big)(t),
\end{eqnarray*}%%----------------------------%
from $\eqref{f2}$, then for $\frac{r_2(t)}{2} \leq r_1(t) \leq 2r_2(t)$
\begin{eqnarray*}%%%%%%%%%%%%%
&&||r_1(t) - r_2(t)||_{L^{2\land \gamma}\big( \big\{\frac{r_2(t)}{2} \leq r_1(t) \leq 2r_2(t) \big\}\big)}^2 \leq c||r_1(t) - r_2(t)||_{L^{2\land \gamma}\big( \big\{\frac{r_2(t)}{2} \leq r_1(t) \leq 2r_2(t) \big\}\big)}^{2\land \gamma} \\
&& \leq c \mathcal{E}\Big((r_1,\bb{U}_1,w_1)\Big|(r_2,\bb{U}_2,w_2) \Big)(t),
\end{eqnarray*}%%----------------------------%%
from $\eqref{f1}$, and finally for $r_1(t) \geq 2r_2(t)$
\begin{eqnarray*}%%%%%%%%%%%%%
&&||r_1(t) - r_2(t)||_{L^{2 \land \gamma}(\{r_1(t) \geq 2r_2(t)\})}^2 =\Big|\Big|\Big(1 - \frac{r_2(t)}{r_1(t)}\Big)r_1(t) \Big|\Big|_{L^{2 \land \gamma}(\{r_1(t) \geq 2 r_2(t)\})}^2 \leq || r_1(t) ||_{L^{2 \land \gamma}(\{r_1(t) \geq 2r_2(t)\})}^2\\
&& \leq C\Big(\int_{\{r_1(t) \geq 2r_2(t)\}}1+r_1^{2 \land \gamma}\Big)^{\frac{2}{2\land \gamma}}\leq C\mathcal{E}\Big((r_1,\bb{U}_1,w_1)\Big|(r_2,\bb{U}_2,w_2) \Big)^{\frac{2}{2 \land \gamma}}(t)\\
&& \leq C\mathcal{E}\Big((r_1,\bb{U}_1,w_1)\Big|(r_2,\bb{U}_2,w_2) \Big)(t),
\end{eqnarray*}%%----------------------------%%
from $\eqref{f2}$, since $\mathcal{E}\Big((r_1,\bb{U}_1,w_1)\Big|(r_2,\bb{U}_2,w_2) \Big)(t) \leq C$, so by combining previous three inequalities, one obtains
\begin{eqnarray}\label{twolandgamma}%%%%%%%%%%%%%
|| r_1(t) - r_2(t)||_{L^{2 \land \gamma}(\Omega)}^{2} \leq C\mathcal{E}\Big((r_1,\bb{U}_1,w_1)\Big|(r_2,\bb{U}_2,w_2) \Big)(t).
\end{eqnarray}%%----------------------------%%

${}$\\
\noindent
\textbf{\textit{Step 3: closing the estimates}}\\

\noindent
Here, we will prove the right-hand side of inequality $\eqref{estimate}$ can be controlled by $\eqref{provethis}$, as follows.\\

\textbf{The 1st and the 2nd term on RHS of }$\eqref{estimate}$:
\begin{eqnarray*}%%%%%%%%%%%%%
	 &&\int_{Q_t} J \mathbb{S}( \nabla_1 \bb{U}_2):(\nabla_1\bb{U}_2 - \nabla_1 \bb{U}_1)  -\int_{Q_t}J_2 \mathbb{S}(\nabla_2\bb{U}_2):(\nabla_2\bb{U}_2 - \nabla_2\bb{U}_1) \\
	 &&=  \int_{Q_t} (J-J_2) \mathbb{S}( \nabla_1 \bb{U}_2):(\nabla_1\bb{U}_2 - \nabla_1 \bb{U}_1) + \int_{Q_t} J_2 \mathbb{S}( (\nabla_1 - \nabla_2) \bb{U}_2):(\nabla_1\bb{U}_2 - \nabla_1 \bb{U}_1) \\
	 &&\quad +\int_{Q_t} J_2 \mathbb{S}( \nabla_2 \bb{U}_2):(\nabla_1-\nabla_2)(\bb{U}_2 - \bb{U}_1) \\
	 &&\leq C\int_0^t||w_1 - w_2||_{L^\infty(\Gamma)} || \nabla \bb{U}_2||_{L^p(\Omega)} ||\nabla A_1||_{L^{\infty^-}(\Omega)}	|| \nabla\bb{U}_1 - \nabla \bb{U}_2||_{L^{2^-}(\Omega)}||\nabla A_1||_{L^{\infty^-}(\Omega)} \\
	 &&\quad+C\int_0^t  || A_1 -  A_2||_{W^{1,\infty^-}(\Omega)} || \nabla \bb{U}_2||_{L^p(\Omega)} ||\nabla A_1||_{L^{\infty^-}(\Omega)}	|| \nabla\bb{U}_1 - \nabla \bb{U}_2||_{L^{2^-}(\Omega)}||\nabla A_1||_{L^{\infty^-}(\Omega)}\\
	&&\quad+C \int_0^t || \nabla \bb{U}_2||_{L^p(\Omega)}|| \nabla A_2||_{L^\infty(\Omega)}||A_1 -  A_2||_{W^{1,\infty^-}(\Omega)} || \nabla\bb{U}_1 - \nabla \bb{U}_2||_{L^{2^-}(\Omega)} \\
	&&\leq C\int_0^t ||\Delta w_1 - \Delta w_2||_{L^2(\Gamma)}|| \bb{U}_1 - \bb{U}_2||_{W^{1,2^-}(\Omega)} \\
	&&\leq C(\delta)\int_0^t\mathcal{E}\Big((r_1,\bb{U}_1,w_1)\Big|(r_2,\bb{U}_2,w_2) \Big) +\delta \int_{Q_t} J |\mathbb{S} (\nabla_1 \bb{U}_1 - \nabla_1 \bb{U}_2|^2,
\end{eqnarray*}%%----------------------------%%
for any $p\in (2,3)$, since $J,J_2 \leq C$ and $\nabla \bb{U}_2 \in L^\infty(0,T; L^p(\Omega))$, where we also used
\begin{eqnarray*}%%%%%%%%%%%%%
 &&|| \bb{U}_1 - \bb{U}_2||_{W^{1,2^-}(\Omega)}^2 \leq C|| \nabla A_1||_{L^{\infty^-}(\Omega)}\int_\Omega J|\mathbb{S}(\nabla_1 \bb{U}_1 - \nabla_1 \bb{U}_2)|^2 \\
 &&\leq C|| \nabla A_1||_{H^{1}(\Omega)}^2\int_\Omega J|\mathbb{S}(\nabla_1 \bb{U}_1 - \nabla_1 \bb{U}_2)|^2 \leq  C\underbrace{|| \Delta w_1||_{L^2(\Omega)}^2}_{\leq C}\int_\Omega J|\mathbb{S}(\nabla_1 \bb{U}_1 - \nabla_1 \bb{U}_2)|^2
\end{eqnarray*}%%----------------------------%%
and
\begin{eqnarray*}%%%%%%%%%%%%%
	||(\nabla_1 - \nabla_2)f||_{L^{p^-}(\Omega)} \leq C ||A_1 - A_2||_{W^{1,\infty^-}(\Omega)}||\nabla f||_{L^p(\Omega)}  \leq C ||\Delta w_1 - \Delta w_2||_{L^{2}(\Gamma)}||\nabla f||_{L^p(\Omega)},
\end{eqnarray*}%%----------------------------%%
for any $f\in W^{1,p}(\Omega)$, by the imbedding of Sobolev spaces.\\

\textbf{The 3rd term on RHS of }$\eqref{estimate}$:
\begin{eqnarray*}%%%%%%%%%%%%%
&&\int_{Q_t}J\frac{1}{r_2}(r_1 - r_2)\big[ \nabla_2 \cdot \mathbb{S}(\nabla_2\bb{U}_2)\big] \cdot (\bb{U}_2 - \bb{U}_1)\\
&&\leq C \int_0^t ||r_1 - r_2||_{L^{2\land \gamma}(\Omega)} \underbrace{|| \nabla_2 \cdot \mathbb{S}(\nabla_2\bb{U}_2) ||_{ L^q(\Omega)}}_{:=h_1 \in L^2(0,T)}|| \bb{U}_2 - \bb{U}_1||_{L^{6^-}(\Omega)} \\
&&\leq C(\delta)\int_0^t h_1^2 ||r_1 - r_2||_{L^{2\land \gamma}(\Omega)}^2 + \delta \int_{Q_t} J |\mathbb{S} (\nabla_1 \bb{U}_1 - \nabla_1 \bb{U}_2)|^2 \\
&&\leq C(\delta) \int_0^t h_1^2 \mathcal{E}\Big((r_1,\bb{U}_1,w_1)\Big|(r_2,\bb{U}_2,w_2) \Big) +\delta \int_{Q_t} J |\mathbb{S} (\nabla_1 \bb{U}_1 - \nabla_1 \bb{U}_2)|^2,
\end{eqnarray*}%%----------------------------%%
for $q>\max \big\{3,\frac{6\gamma}{5\gamma-6} \big\}$, by using $\eqref{twolandgamma}$ and the following inequalities
\begin{eqnarray*}%%%%%%%%%%%%%
&&|| \nabla_2 \cdot \mathbb{S}(\nabla_2\bb{U}_2) ||_{ L^q(\Omega)}\leq C\Big( || \nabla^2 \bb{U}_2||_{L^{q}(\Omega)}||(\nabla w_2)^2||_{L^{\infty}(\Gamma)}+ || \nabla \bb{U}_2||_{L^{\infty}(\Omega)}||\Delta w_2||_{L^{\infty}(\Gamma)} \Big), \\[2mm]
&&|| \bb{U}_1 - \bb{U}_2||_{L^{6^-}(\Omega)}^2 \leq C || \bb{U}_1 - \bb{U}_2||_{W^{1,2^-}(\Omega)}^2 \leq C\underbrace{|| \nabla A_1||_{L^{\infty^-}(\Omega)}^2}_{\leq C}\int_\Omega J|\mathbb{S}(\nabla_1 \bb{U}_1 - \nabla_1 \bb{U}_2)|^2.
\end{eqnarray*}%%----------------------------%%

\textbf{The 4th term on RHS of }$\eqref{estimate}$:
\begin{eqnarray*}%%%%%%%%%%%%%
	&& \int_{Q_t}(J-J_2)\nabla_2 \cdot \mathbb{S}(\nabla_2\bb{U}_2)(\bb{U}_2 - \bb{U}_1) \leq \int_0^t ||w_1 - w_2||_{L^\infty(\Gamma)} h_1 || \bb{U}_1 - \bb{U}_2||_{L^{6^-}(\Omega)}\\
	 &&\leq C(\delta) \int_0^t h_1^2 \mathcal{E}\Big((r_1,\bb{U}_1,w_1)\Big|(r_2,\bb{U}_2,w_2) \Big) +\delta \int_{Q_t} J |\mathbb{S} (\nabla_1 \bb{U}_1 - \nabla_1 \bb{U}_2)|^2 .
\end{eqnarray*}%%----------------------------%%

\textbf{The 5th term on RHS of }$\eqref{estimate}$ - \textbf{case} $\alpha = 0, \gamma>12/7$:
\begin{eqnarray*}%%%%%%%%%%%%%
	 && \int_{Q_t} J(r_2 - r_1 )\Big[ \bb{U}_2 \cdot (\nabla_1 - \nabla_2)(r_2^{\gamma-1}) - (\bb{w}_1 - \bb{w}_2)\cdot\nabla_1(r_2^{\gamma-1}) - \bb{w}_2\cdot(\nabla_1 - \nabla_2)(r_2^{\gamma-1}) \Big] \\
	&& \leq C \int_0^t ||r_1-r_2||_{L^{2\land\gamma}(\Omega)} \Big[||\Delta w_1- \Delta w_2||_{L^{2}(\Gamma)}\underbrace{||\bb{U}_2 \cdot\nabla (r_2^{\gamma-1}) ||_{L^q(\Omega)}}_{:=h_2 \in L^2(0,T)} + ||\partial_t w_1 - \partial w_2||_{L^{4^-}(\Gamma)}\underbrace{|| \nabla (r_2^{\gamma-1} )||_{L^q(\Omega)}}_{:=h_3\in L^2(0,T)}\\
&&\quad \quad \quad\quad \quad \quad\quad \quad \quad~+||\Delta w_1 - \Delta w_2||_{L^{2}(\Gamma)}\underbrace{||\bb{w}_2 \cdot\nabla (r_2^{\gamma-1}) ||_{L^q(\Omega)}}_{:=h_4\in L^2(0,T)} \Big]  \\
	&&\leq C(\delta)\int_0^t (h_2^2+h_3^2+h_4^2) ||r_1-r_2||_{L^{2\land\gamma}(\Omega)}^2+ C\int_0^t ||\Delta w_1 - \Delta w_2||_{L^2(\Gamma)}^2 +\delta \int_0^t ||\partial_t w_1 - \partial w_2||_{H^{({\frac{1}{2}})^-}(\Gamma)}^2 \\
	&&\leq C(\delta)\int_0^t(h_2^2 +h_3^2+ h_4^2+1)\mathcal{E}\Big((r_1,\bb{U}_1,w_1)|(r_2,\bb{U} _2,w_2)\Big)+\delta \int_{Q_t} J |\mathbb{S} (\nabla_1 \bb{U}_1 - \nabla_1 \bb{U}_2)|^2 ,
\end{eqnarray*}%%----------------------------%%
for any $q>\max \{4, \frac{4\gamma}{3\gamma-4} \}$, where we used 
\begin{eqnarray}%%%%%%%%%%%%%
&&||\partial_t w_1 - \partial_t w_2||_{L^{4^-}(\Gamma)}^2 \leq C||\partial_t w_1 - \partial_t w_2||_{W^{(\frac{1}{2})^-, 2^-}(\Gamma)}^2 \leq C || \bb{U}_1 - \bb{U}_2||_{W^{1,2^-}(\Omega)}^2\nonumber\\[2mm]
&& \leq C\int_{\Omega} J |\mathbb{S} (\nabla_1 \bb{U}_1 - \nabla_1 \bb{U}_2)|^2,\label{tracedis}
\end{eqnarray}%%----------------------------%%
by imbedding and trace inequality.

\textbf{The 5th term on RHS of }$\eqref{estimate}$ - \textbf{case} $\alpha>0,\gamma>3/2$:
\begin{eqnarray*}%%%%%%%%%%%%%
	&& \int_{Q_t} J(r_2 - r_1 )\Big[ \bb{U}_2 \cdot (\nabla_1 - \nabla_2)(r_2^{\gamma-1}) - (\bb{w}_1 - \bb{w}_2)\cdot\nabla_1(r_2^{\gamma-1}) - \bb{w}_2\cdot(\nabla_1 - \nabla_2)(r_2^{\gamma-1}) \Big]\\
	&&\leq C \int_0^t ||r_1-r_2||_{L^{2\land\gamma}(\Omega)} \Big[||\Delta w_1 -  \Delta w_2||_{L^{2}(\Gamma)}\underbrace{||\bb{U}_2 \cdot\nabla (r_2^{\gamma-1}) ||_{L^q(\Omega)}}_{:=h_2\in L^2(0,T)} + ||\partial_t w_1 - \partial w_2||_{L^{\infty^-}(\Gamma)}\underbrace{|| \nabla (r_2^{\gamma-1}) ||_{L^q(\Omega)}}_{:=h_3\in L^2(0,T)}\\
	&&\quad \quad \quad\quad \quad \quad\quad \quad \quad~+||\Delta w_1 - \Delta w_2||_{L^{2}(\Gamma)}\underbrace{||\bb{w}_2 \cdot \nabla (r_2^{\gamma-1}) ||_{L^q(\Omega)}}_{:=h_4\in L^2(0,T)} \Big]  \\
	&&\leq C(\delta)\int_0^t(h_2^2 +h_3^2+ h_4^2+1)\mathcal{E}\Big((r_1,\bb{U}_1,w_1)\Big|(r_2,\bb{U}_2,w_2) \Big)+\delta\alpha \int_0^t ||\partial_t \nabla w_1 - \partial_t \nabla w_2||_{L^2(\Gamma)}^2,
\end{eqnarray*}%%----------------------------%%
for any $q>3$.\\

\textbf{The 6th term on RHS of }$\eqref{estimate}$:
\begin{eqnarray*}%%%%%%%%%%%%%
	 &&\int_{Q_t} Jr_1 \big[ \bb{U}_2\cdot (\nabla_1 - \nabla_2)\bb{U}_2\big]\cdot(\bb{U}_2 - \bb{U}_1) \\
	 &&\leq C \int_0^t||r_1||_{L^{\gamma}(\Omega)} ||\bb{U}_2||_{L^\infty(\Omega)} ||\Delta w_1 - \Delta w_2||_{L^{2}(\Gamma)} \underbrace{|| \nabla \bb{U}_2||_{L^\infty(\Omega)}}_{=h_5 \in L^2(0,T)} || \bb{U}_1 - \bb{U}_2||_{ L^{6^-}(\Omega)} \\
	&&\leq C(\delta) \int_0^t h_5^2 \mathcal{E}\Big((r_1,\bb{U}_1,w_1)\Big|(r_2,\bb{U}_2,w_2) \Big) + \delta \int_{Q_t} J |\mathbb{S} (\nabla_1 \bb{U}_1 - \nabla_1 \bb{U}_2)|^2 .
\end{eqnarray*}%%----------------------------%%

\textbf{The 7th term on RHS of }$\eqref{estimate}$:
\begin{eqnarray*}%%%%%%%%%%%%%
	&&\int_{Q_t} Jr_1 \big[(\bb{U}_1 - \bb{U}_2) \cdot \nabla_1 \bb{U}_2\big]\cdot(\bb{U}_2 - \bb{U}_1) \\
	&&\leq \int_{0}^t \Big[\underbrace{|| \nabla w_1||_{L^\infty(\Gamma)}|| \nabla \bb{U}_2||_{L^\infty(\Omega)}}_{:=h_6\in L^1(0,T)}\int_{\Omega} Jr_1 |\bb{U}_2 - \bb{U}_1|^2\Big] \leq C\int_0^t h_6 \mathcal{E}\Big((r_1,\bb{U}_1,w_1)\Big|(r_2,\bb{U}_2,w_2) \Big),
\end{eqnarray*}%%----------------------------%%
since $\nabla w_1 \in L^2(0,T;L^\infty(\Gamma))$, by Theorem $\ref{main2}(1)$.\\

\textbf{The 8th term on RHS of }$\eqref{estimate}$:
\begin{eqnarray*}%%%%%%%%%%%%%	
	&&\int_{Q_t}J r_1(\nabla_1 - \nabla_2)(r_2^{\gamma-1})\cdot(\bb{U}_2 - \bb{U}_1) \leq C\int_0^t \underbrace{||r_1 ||_{L^\gamma(\Omega)}||\nabla r_2 ||_{L^p(\Omega)}}_{:=h_7\in L^2(0,T)} || \Delta w_1 - \Delta w_2||_{L^2(\Gamma)}|| \bb{U}_1 - \bb{U}_2||_{L^{6^-}(\Omega)}  \\
	&&\leq C(\delta) \int_0^t h_7^2 \mathcal{E}\Big((r_1,\bb{U}_1,w_1)\Big|(r_2,\bb{U}_2,w_2) \Big) +\delta \int_{Q_t} J |\mathbb{S} (\nabla_1 \bb{U}_1 - \nabla_1 \bb{U}_2)|^2 ,
\end{eqnarray*}%%----------------------------%%
for any $p>\frac{6\gamma}{5\gamma-6}$. \\

\textbf{The 9th term on RHS of }$\eqref{estimate}$:
\begin{eqnarray*}%%%%%%%%%%%%%
	&&\int_{Q_t} J (r_2 - r_1) r_2^{\gamma-1}(\nabla_1 - \nabla_2) \cdot \bb{U}_2 \leq C\int_0^t \underbrace{||\nabla \bb{U}_2||_{L^\infty(\Omega)}}_{:=h_8 \in L^2(0,T)}|| r_1 - r_2||_{L^{2\land \gamma}(\Omega)} || \Delta w_1 - \Delta  w_2||_{L^{2}(\Gamma)}  \\
	&&\leq C  \int_0^t h_8 \mathcal{E}\Big((r_1,\bb{U}_1,w_1)\Big|(r_2,\bb{U}_2,w_2) \Big).
\end{eqnarray*}%%----------------------------%%

\textbf{The 10th term on RHS of }$\eqref{estimate}$:
\begin{eqnarray*}%%%%%%%%%%%%%
	&&\int_{Q_t} J\big[ r_2^\gamma - \gamma r_2^{\gamma-1}(r_2 - r_1) - r_1^\gamma\big](\nabla_1 \cdot \bb{U}_2)\\
	&&\leq \int_0^t \Big[\underbrace{|| \nabla w_1||_{L^\infty(\Gamma)}|| \nabla \bb{U}_2||_{L^\infty(\Omega)}}_{:=h_9\in L^1(0,T)}\int_{\Omega}J\big( r_2^\gamma - \gamma r_2^{\gamma-1}(r_2 - r_1) - r_1^\gamma\big) \Big] \leq  C \int_0^t h_9 \mathcal{E}\Big((r_1,\bb{U}_1,w_1)\Big|(r_2,\bb{U}_2,w_2) \Big).
\end{eqnarray*}%%----------------------------%%

\textbf{The 11th term on RHS of }$\eqref{estimate}$ - \textbf{case} $\alpha = 0$, $\gamma \in (\frac{12}{7},2]$:\\

Assuming additionally that $\nabla \bb{U}_2 \in L^\infty(Q_T)$, one can bound
\begin{eqnarray*}%%%%%%%%%%%%%
&&\int_{Q_t} Jr_1\Big[\bb{w}_1 \cdot \nabla_1 \bb{U}_2 -\bb{w}_2\cdot  \nabla_2 \bb{U}_2\Big]\cdot (\bb{U}_1 - \bb{U}_2) \\
&&     \int_{Q_t} J(r_1-r_2)\Big[\bb{w}_1 \cdot \nabla_1 \bb{U}_2 -\bb{w}_2\cdot  \nabla_2 \bb{U}_2\Big]\cdot (\bb{U}_1 - \bb{U}_2)    +\int_{Q_t} Jr_2\Big[\bb{w}_1 \cdot \nabla_1 \bb{U}_2 -\bb{w}_2\cdot  \nabla_2 \bb{U}_2\Big]\cdot (\bb{U}_1 - \bb{U}_2)           \\
&&\leq C \int_0^t ||r_1-r_2||_{L^{ \gamma}(\Omega)}||\nabla \bb{U}_2 ||_{L^\infty(\Omega)}\underbrace{\Big[  || \partial_t w_1 \nabla w_1||_{L^{4^-}(\Gamma)}+||\partial_t w_2 \nabla w_2||_{L^\infty(\Gamma )}  \Big]}_{:=h_{10}\in L^2(0,T)}||\bb{U}_1 - \bb{U}_2||_{L^{6^-}(\Omega)}\\
&&\quad+C \int_0^t ||r_2||_{L^{\infty}(\Omega)}||\nabla \bb{U}_2 ||_{L^\infty(\Omega)}\Big[|| \nabla w_1 ||_{L^{\infty^-}(\Gamma)}  ||\partial_t w_1-\partial_t w_2||_{L^{2}(\Gamma)}\\
&&\quad \quad\quad  \quad\quad  \quad\quad  \quad\quad  \quad\quad  \quad\quad  \quad +||\partial_t w_2||_{L^\infty(\Gamma)}|| \Delta w_1 - \Delta w_2||_{L^{2}(\Gamma)} \Big]||\bb{U}_1 - \bb{U}_2||_{L^{6^-}(\Omega)} \\
&&\leq C(\delta)\int_0^t \Big[h_{10}^2||r_1-r_2||_{L^{ \gamma}(\Omega)}^2+||\partial_t w_1-\partial_t w_2||_{L^2(\Gamma)}^2+||\Delta w_1-\Delta w_2||_{L^2(\Gamma)}^2\Big]+\delta \int_{Q_t} J |\mathbb{S} (\nabla_1 \bb{U}_1 - \nabla_1 \bb{U}_2)|^2 \\
&&\leq C(\delta) \int_0^t (h_{10}^2+1) \mathcal{E}\Big((r_1,\bb{U}_1,w_1)\Big|(r_2,\bb{U}_2,w_2) \Big) +\delta \int_{Q_t} J |\mathbb{S} (\nabla_1 \bb{U}_1 - \nabla_1 \bb{U}_2)|^2. \\
\end{eqnarray*}%%----------------------------%%

\textbf{The 11th term on RHS of }$\eqref{estimate}$ - \textbf{case} $\alpha = 0$, $\gamma>2$:
\begin{eqnarray*}%%%%%%%%%%%%%
	&&\int_{Q_t} Jr_1\Big[\bb{w}_1 \cdot \nabla_1 \bb{U}_2 -\bb{w}_2\cdot  \nabla_2 \bb{U}_2\Big]\cdot (\bb{U}_1 - \bb{U}_2) \\
	&&\leq C \int_0^t ||\sqrt{r_1}||_{L^{2\gamma}(\Omega)}\underbrace{||\nabla \bb{U}_2 ||_{L^\infty(\Omega)}}_{:=h_{11}\in L^2(0,T)}\Big[|| \nabla w_1 ||_{L^{\infty^-}(\Gamma)}  ||\partial_t w_1-\partial_t w_2||_{L^{4^-}(\Gamma)}\\
	&&\quad \quad\quad  \quad\quad  \quad\quad  \quad\quad  \quad\quad  \quad\quad  \quad +||\partial_t w_2||_{L^\infty(\Gamma)}||\Delta  w_1 -  \Delta w_2||_{L^{2}(\Gamma)} \Big]||\sqrt{r_1}(\bb{U}_1 - \bb{U}_2)||_{L^{2}(\Omega)} \\
	&&\leq C\int_0^t ||\Delta w_1-\Delta w_2||_{L^2(\Gamma)}^2+C(\delta)\int_0^th_{11}^2 \int_{\Omega}J r_1|\bb{U}_1 - \bb{U}_2|^2+\delta \int_{Q_t} J |\mathbb{S} (\nabla_1 \bb{U}_1 - \nabla_1 \bb{U}_2)|^2\\[2mm]
	&&\leq C(\delta) \int_0^t (h_{11}^2+1)\mathcal{E}\Big((r_1,\bb{U}_1,w_1)\Big|(r_2,\bb{U}_2,w_2) \Big)+\delta \int_{Q_t} J |\mathbb{S} (\nabla_1 \bb{U}_1 - \nabla_1 \bb{U}_2)|^2,
\end{eqnarray*}%%----------------------------%%
by $\eqref{tracedis}$.\\

\textbf{The 11th term on RHS of }$\eqref{estimate}$ - \textbf{case} $\alpha > 0$, $\gamma>3/2$:
\begin{eqnarray*}%%%%%%%%%%%%%
	&&\int_{Q_t} Jr_1\Big[\bb{w}_1 \cdot \nabla_1 \bb{U}_2 -\bb{w}_2\cdot  \nabla_2 \bb{U}_2\Big]\cdot (\bb{U}_1 - \bb{U}_2) \\
	&&\leq C \int_0^t ||\sqrt{r_1}||_{L^{3}(\Omega)}\underbrace{||\nabla \bb{U}_2 ||_{L^\infty(\Omega)}}_{=h_{12}\in L^2(0,T)}\Big[|| \nabla w_1 ||_{L^{\infty^-}(\Gamma)}  ||\partial_t w_1-\partial_t w_2||_{L^{\infty^-}(\Gamma)}\\
	&&\quad \quad\quad  \quad\quad  \quad\quad  \quad\quad  \quad\quad  \quad\quad  \quad +||\partial_t w_2||_{L^\infty(\Gamma)}||\Delta w_1 -  \Delta w_2||_{L^{2}(\Gamma)} \Big]||\sqrt{r_1}(\bb{U}_1 - \bb{U}_2)||_{L^{2}(\Omega)} \\
	&&\leq C\int_0^t ||\Delta w_1-\Delta w_2||_{L^2(\Gamma)}^2+C(\delta)\int_0^t h_{12}^2 \int_{\Omega}J r_1|\bb{U}_1 - \bb{U}_2|^2+ \delta \alpha \int_0^t ||\partial_t\nabla  w_1-\partial_t \nabla  w_2||_{L^2(\Gamma)}^2\\
	&&\leq C \int_0^t (h_{12}^2+1)\mathcal{E}\Big((r_1,\bb{U}_1,w_1)\Big|(r_2,\bb{U}_2,w_2) \Big)+ \delta \alpha \int_0^t ||\partial_t\nabla  w_1-\partial_t \nabla  w_2||_{L^2(\Gamma)}^2.
\end{eqnarray*}%%----------------------------%%

Combining previous inequalities, the proof is finished.

\section{Some nonlinear plate models}
Here, we study the interaction between compressible viscous fluids and nonlinear plates. The goal is to prove that the main results of this paper can be extended to such interaction problems. 
We note that the existence of weak solutions for the following interaction problems were constructed in \cite{trwa3}, while the existence of strong solutions seem to be non-existent in the literature. 

\subsection{A semilinear elastic plate}\label{semi}
Here, the equation $\eqref{structureeqs}$ is replaced by
\begin{eqnarray*}%%%%%%%%%%%%%
\partial_t^2 w+\Delta^2 w - \alpha \partial_t \Delta w + \mathcal{F}(w)=-S^w \bb{f}_{fl} \cdot \mathbf{e_3},
\end{eqnarray*}%%----------------------------%%
where $\mathcal{F}$ satisfies the following assumptions:
\begin{enumerate}
\item[(A1)] Mapping $\mathcal{F}$ is locally Lipschitz from $H_0^{2-\epsilon}(\Gamma)$ into $H^{-2}(\Gamma)$ for some $\epsilon >0$, i.e.
\begin{eqnarray*}
||\mathcal{F}(w_1)-\mathcal{F}(w_2)||_{H^{-2}(\Gamma)} \leq
C_R || w_1 -w_2 ||_{H^{2-\epsilon}(\Gamma)},
\end{eqnarray*}
for a constant $C_R>0$, for any $||w_i||_{H^{2}(\Gamma)} \leq R$ ($i=1,2$). 

\item[(A2)] Mapping $\mathcal{F}$ is locally Lipschitz from $H_0^{2}(\Gamma)$ into $H^{-a}(\Gamma)$ for some $0\leq a<1/2$, i.e.
\begin{eqnarray}\label{assumption1}
||\mathcal{F}(w_1)-\mathcal{F}(w_2)||_{H^{-a}(\Gamma)} \leq
C_R || w_1 -w_2 ||_{H^{2}(\Gamma)},
\end{eqnarray}
for a constant $C_R>0$, for any $||w_i||_{H^2(\Gamma)} \leq R$ ($i=1,2$). 

\item[(A3)]  $\mathcal{F}(w)$ has a potential
in $H_0^2(\Gamma)$, i.e. there exists a Fr\'{e}chet
differentiable functional $\Pi(w)$ on $H_0^2(\Gamma)$ such that
$\Pi'(w)=\mathcal{F}(w)$, and there are $0<\kappa <1/2$ and  $C^* \geq 0$, such that the following inequality holds,
\begin{eqnarray*}
 \kappa || \Delta
w||_{L^2(\Gamma)}^2+\Pi(w)+ C^*\geq 0, ~~ \text{for all } w \in H_0^2(\Gamma).
\end{eqnarray*}
Moreover, $\Pi$ is bounded on bounded sets in $H_0^2(\Gamma)$.
\end{enumerate}

\begin{rem}
(1) This plate model is a generalization of Kirchhoff, von Karman and Berger plates which are presented in Appendix. \\
(2)	In \cite{trwa3}, the condition $(A2)$ was not necessary for the existence proof and was thus removed. However, in the proof that follows, it is needed. 
\end{rem}

Here, we give a sketch of the proofs of theorems $\ref{main2}$ and $\ref{main1}$ for this plate model. Let $(\rho_1,\bb{u}_1,w_1)$ be a finite energy weak solution. In this case, the relative entropy stays the same, while in the relative energy inequality takes the form
\begin{eqnarray}%%%%%%%%%%%%%
&&\mathcal{E}\Big((\rho_1,\bb{u}_1,w_1)\Big|(\rho_2,\bb{u}_2,w_2)\Big)(t)  + \int_{0}^t\int_{\Omega^{w_1}(t)}  \mathbb{S} (\nabla \bb{u}_1 - \nabla \bb{u}_2):(\nabla\bb{u}_1 - \nabla \bb{u}_2) + \alpha\int_0^t\int_{\Gamma}  |\partial_t \nabla w_1 - \partial_t \nabla w_2|^2 \nonumber\\
&&\leq \mathcal{E}\Big((\rho_1,\bb{u}_1,w_1)\Big|(\rho_2,\bb{u}_2,w_2)\Big)(0) +\int_0^t\mathcal{R}\Big( \rho_1,\bb{u}_1,w_1, \rho_2,\bb{u}_2,w_2\Big) -\int_0^t \langle \mathcal{F}(w_1), \partial_t w_1 - \partial_t w_2 \rangle \nonumber
\end{eqnarray}%%----------------------------%%
It is easy to prove that any finite energy weak solution satisfies this inequality. Next, to obtain the regularity estimates in  Theorem $\ref{main2}(1)$, in $\eqref{aaaa}$ additional term $-\int_{\Gamma_T} \langle \mathcal{F}(w), D_h^s w \rangle$ appears, which can be controlled as
\begin{eqnarray*}%%%%%%%%%%%%%
&&-\int_0^t \langle \mathcal{F}(w), D_{-h}^s D_h^s w \rangle \leq \int_0^t \Big[|| \mathcal{F}(w) - \mathcal{F}(0)||_{H^{-a}(\Gamma)}+||\mathcal{F}(0)||_{H^{-a}(\Gamma)} \Big] || D_{-h}^s D_h^s w||_{H^{a}(\Gamma)} \\
&&\leq \int_0^t C_R(||w||_{H^2(\Gamma)}+1)||  w||_{H^2(\Gamma)}\leq C
\end{eqnarray*}%%----------------------------%%
by $(A2)$ and for any $s<\frac{1}{2}$, where we have used the fact that the constant $C_R$ in $\eqref{assumption1}$ is uniform, since $w_1,w_2$ are uniformly bounded in $L^\infty(0,T;H_0^2(\Gamma))$.

Finally, to obtain the Theorem $\ref{main1}$ in this case, in $\eqref{estimate}$ an additional term appears $-\int_{\Gamma_T} \langle \mathcal{F}(w_1) -\mathcal{F}(w_2),\partial_t w_1 - \partial_t w_2\rangle$ on the right-hand side, which can be controlled as
\begin{eqnarray*}%%%%%%%%%%%%%
&&-\int_{0}^t \langle \mathcal{F}(w_1) -\mathcal{F}(w_2)  , \partial_t w_1 - \partial_t w_2\rangle \leq \int_0^t||\mathcal{F}(w_1) -\mathcal{F}(w_2) ||_{H^{-a}(\Gamma)} || \partial_t w_1 - \partial_t w_2||_{H^{a}(\Gamma)} \\
&&\leq C \int_0^t ||w_1 -w_2 ||_{H^{2}(\Gamma)} || \partial_t w_1 - \partial_t w_2||_{H^{a}(\Gamma)}\\
&& \leq C(\delta)\int_0^t \mathcal{E}\Big(\big[(\rho_1,\bb{u}_1,w_1)|(\rho_2,\bb{u}_2,w_2)\big] \Big) +\delta \int_{Q_t} J |\mathbb{S} (\nabla_1 \bb{U}_1 - \nabla_1 \bb{U}_2)|^2,
\end{eqnarray*}%%----------------------------%%
by $\eqref{tracedis}$. If $\alpha>0$, then we can choose $a\leq 1$ and bound
\begin{eqnarray*}%%%%%%%%%%%%%
	&&-\int_{\Gamma_T} \langle \mathcal{F}(w_1) -\mathcal{F}(w_2)  , \partial_t w_1 - \partial_t w_2\rangle \leq \int_0^t||\mathcal{F}(w_1) -\mathcal{F}(w_2) ||_{H^{-a}(\Gamma)} || \partial_t w_1 - \partial_t w_2||_{H^{a}(\Gamma)} \\
	&&\leq C \int_0^t ||w_1 -w_2 ||_{H^{2}(\Gamma)} || \partial_t w_1 - \partial_t w_2||_{H^{1}(\Gamma)}\\
	&& \leq C(\delta)\int_0^t \mathcal{E}\Big(\big[(\rho_1,\bb{u}_1,w_1)|(\rho_2,\bb{u}_2,w_2)\big] \Big) + \alpha\delta\int_0^t || \partial_t w_1 - \partial_t w_2||_{H^{1}(\Gamma)}.
\end{eqnarray*}%%----------------------------%%

\subsection{A thermoelastic semilinear plate}
Here, the equation $\eqref{structureeqs}$ is replaced by a system\footnote{Note that in this case, there is an additional equation in the weak formulation.}
\begin{eqnarray*}%%%%%%%%%%%%%
	\partial_t^2 w+\Delta^2 w - \alpha \partial_t \Delta w + \mathcal{F}(w)&=&S^w \bb{f}_{fl} \cdot \mathbf{e_3},\\
	\partial_t \theta - \Delta\theta -\partial_t \Delta w &=& 0,
\end{eqnarray*}%%----------------------------%%
where $\theta:\Gamma \to \mathbb{R}$ is the temperature of the plate and $\mathcal{F}$ satisfies the assumptions $(A1)-(A3)$. Here, the relative entropy takes the form
\begin{align}%%%%%%%%%%%%%
\mathcal{E}\Big((\rho_1,\bb{u}_1,w_1)\Big|(\rho_2,\bb{u}_2,w_2)\Big)(t)  &= \frac{1}{2} \int_{\Omega^{w_1}(t)} (\rho_1 |\bb{u}_1 - \bb{u}_2|^2)(t) + \frac{1}{\gamma-1}\int_{\Omega^{w_1}(t)} \big( \rho_1^\gamma - \gamma \rho_2^{\gamma-1}(\rho_1 - \rho_2) - \rho_2^{\gamma} \big)(t) \nonumber \\
&\quad+ \frac{1}{2}\int_\Gamma |\partial_t w_1 - \partial_t w_2|^2(t) + \frac{1}{2}\int_\Gamma |\Delta w_1 - \Delta w_2|^2(t)+\frac{1}{2}\int_\Gamma |\theta_1 - \theta_2|^2(t) \quad\quad\quad\quad \nonumber
\end{align}%%----------------------------%%
while the relative energy inequality is replaced by
\begin{eqnarray}%%%%%%%%%%%%%
&&\mathcal{E}\Big((\rho_1,\bb{u}_1,w_1)\Big|(\rho_2,\bb{u}_2,w_2)\Big)(t)  + \int_{0}^t\int_{\Omega^{w_1}(t)} J \mathbb{S} (\nabla \bb{u}_1 - \nabla \bb{u}_2):(\nabla\bb{u}_1 - \nabla \bb{u}_2) \nonumber\\
&&\quad + \alpha\int_0^t\int_{\Gamma}  |\partial_t \nabla w_1 - \partial_t \nabla w_2|^2+ \int_0^t\int_{\Gamma}  | \nabla \theta_1 -  \nabla \theta_2|^2\nonumber\\[2mm]
&&\leq \mathcal{E}\Big(\big[(\rho_1,\bb{u}_1,w_1)|(\rho_2,\bb{u}_2,w_2)\big] \Big)(0) +\int_0^t\mathcal{R}\Big( \rho_1,\bb{u}_1,w_1, \rho_2,\bb{u}_2,w_2\Big)-\int_0^t \langle \mathcal{F}(w_1), \partial_t w_1 - \partial_t w_2\rangle \nonumber\\
&&\quad -\int_0^t \int_{\Gamma}(\theta_2 - \theta_1)\partial_t\theta_2 -\int_0^t \int_{\Gamma} (\theta_2 - \theta_1) \Delta \theta_2.  \nonumber
\end{eqnarray}%%----------------------------
The proofs of theorems $\ref{main2}$ and $\ref{main1}$ in this case can be carried out in the same way as for the semilinear plate case so we omit them here.

\subsection{A thermoelastic quasilinear plate}
Here, the equation $\eqref{structureeqs}$ is replaced by a system
\begin{eqnarray*}%%%%%%%%%%%%%
	\partial_t^2 w+\Delta^2 w - \alpha \partial_t \Delta w + \Delta(\Delta w)^3&=&S^w \bb{f}_{fl} \cdot \mathbf{e_3},\\
	\partial_t \theta - \Delta\theta -\partial_t \Delta w &=& 0.
\end{eqnarray*}%%----------------------------%%
\begin{rem}
	Such a thermoelastic plate model was first studied in \cite{thefirstpaper} and later in \cite{lasiecka}. This high order nonlinearity comes from a thermoelastic plate model in which a nonlinear coupling between the elastic, magnetic and thermoelastic fields is considered.
\end{rem}

In this case, the relative entropy takes the form
\begin{eqnarray}%%%%%%%%%%%%%
    \mathcal{E}\Big((\rho_1,\bb{u}_1,w_1)\Big|(\rho_2,\bb{u}_2,w_2)\Big)(t)  &=& \frac{1}{2} \int_{\Omega^{w_1}(t)} (\rho_1 |\bb{u}_1 - \bb{u}_2|^2)(t) + \frac{1}{\gamma-1}\int_{\Omega^{w_1}(t)} \big( \rho_1^\gamma - \gamma \rho_2^{\gamma-1}(\rho_1 - \rho_2) - \rho_2^{\gamma} \big)(t) \nonumber \\
	&+& \frac{1}{2}\int_\Gamma |\partial_t w_1 - \partial_t w_2|^2(t) + \frac{1}{2}\int_\Gamma |\Delta w_1 - \Delta w_2|^2(t)+\frac{1}{2}\int_\Gamma |\theta_1 - \theta_2|^2(t)\nonumber \\
	&+&\int_\Gamma \Big[(\Delta w_1)^3 - (\Delta w_2)^3\Big] (\Delta w_1- \Delta w_2), \nonumber
\end{eqnarray}%%----------------------------%%
while the relative energy inequality is replaced by
\begin{eqnarray}%%%%%%%%%%%%%
	&&\mathcal{E}\Big((\rho_1,\bb{u}_1,w_1)\Big|(\rho_2,\bb{u}_2,w_2)\Big)(t) + \int_{0}^t\int_{\Omega^{w_1}(t)}  J\mathbb{S} (\nabla \bb{u}_1 - \nabla \bb{u}_2):(\nabla\bb{u}_1 - \nabla \bb{u}_2) \nonumber\\
	&&\quad + \alpha\int_0^t\int_{\Gamma}  |\partial_t \nabla w_1 - \partial_t \nabla w_2|^2+ \int_0^t\int_{\Gamma}  | \nabla \theta_1 -  \nabla \theta_2|^2\nonumber\\[2mm]
	&&\leq \mathcal{E}\Big(\big[(\rho_1,\bb{u}_1,w_1)|(\rho_2,\bb{u}_2,w_2)\big] \Big)(0) +\int_0^t\mathcal{R}\Big( \rho_1,\bb{u}_1,w_1, \rho_2,\bb{u}_2,w_2\Big)-\int_0^t \int_\Gamma \mathcal(\Delta w_2)^3 (\Delta w_1 - \Delta w_2) \nonumber\\
	&&\quad -\int_0^t \int_{\Gamma}(\theta_2 - \theta_1)\partial_t\theta_2 -\int_0^t \int_{\Gamma} (\theta_2 - \theta_1) \Delta \theta_2. \nonumber
\end{eqnarray}%%----------------------------%%
In this case, in the proof of Theorem $\ref{main1}(1)$, when we choose $(\bb{q},\psi) = (\bb{R}[D_{-h}^s D_h^{s}w], D_{-h}^s D_h^{s} w)$ in the coupled momentum equation $\eqref{weaksolmom}$, the nonlinear term can be expressed as
\begin{eqnarray*}%%%%%%%%%%%%%
	\int_{\Gamma_T} (\Delta w)^3 D_{-h}^s D_h^{s}[\Delta w] = \int_{\Gamma_T} D_h^s(\Delta w)^3 D_h^{s}\Delta w \geq \frac{1}{4} \int_{\Gamma_T} |D_h^{\frac{s}{2}} \Delta w|^4.
\end{eqnarray*}%%----------------------------%%
This then gives us $\Delta w \in L^4(0,T; W^{\frac{s}{2},4}(\Gamma))\cap L^2(0,T; H^s(\Gamma))$, where $s$ satisfies the same conditions as in Theorem $\ref{main2}(1)$.

\vspace{.1in}
\noindent{\bf Acknowledgments:} This research was partially carried out during the author's PhD studies at Shanghai Jiao Tong University, where it was supported by the National Natural Science Foundation of China (NNSFC) under Grant No. 11631008. The author is immensely grateful to the referees for their very detailed reading of the manuscript and for many useful comments and suggestions.

\section*{Appendix - the nonlinear plate models}
The assumptions $(A1)-(A3)$ given in section $\ref{semi}$ are satisfied by the following models:\\

\noindent
\textbf{The Kirchhoff model}: Here the nonlinear elastic force takes the form of the Nemytskii operator:
\begin{eqnarray*}%%%%%%%%%%%%%
\mathcal{F}(w) = - \nu  \cdot \text{div} \big[ |\nabla w|^q  \nabla w - \mu |\nabla w|^r \nabla w \big] + f(w) -h(x),
\end{eqnarray*}%%----------------------------%%
where $\nu \geq 0$, $q>r \geq 0$, and $\mu \in \mathbb{R}$ are parameters, $h\in L^2(\Gamma)$ and
\begin{eqnarray*}
f \in \text{Lip}_{loc}(\mathbb{R}) \quad \text{satisfies} \quad \lim\limits_{|s|\to \infty}\inf f(s)s^{-1} > -\lambda_1^2, 
\end{eqnarray*}%%----------------------------%%
where $\lambda_1$ is the first eigenvalue for the Laplacian with the Dirichlet boundary condition. The corresponding potential is
\begin{eqnarray*}%%%%%%%%%%%%%
\Pi(w) = \int_\Gamma \Phi(w) + \frac{\nu}{q+2} \int_{\Gamma} |\nabla w|^{q+2} - \frac{\nu \mu}{r+2} \int_\Gamma |\nabla w|^{r+2} - \int_\Gamma w h,
\end{eqnarray*}%%----------------------------%%
where $\Phi(s) = \int_0^s f(x) dx$. This model is used for the plates with small displacements and it is important for pre-buckling analysis. See \cite{platesproofs} for more details about this model. \\ \\
\textbf{The von K\'{a}rm\'{a}n model}: Here, the nonlinear elastic force takes the form
\begin{eqnarray*}%%%%%%%%%%%%%
\mathcal{F}(w) = -[w, v(w)+F_0] - h
\end{eqnarray*}%%----------------------------%%$, 
where $F_0 \in H^4(\Gamma)$ and $h \in L^2(\Gamma)$ are given functions, and the von Karman bracket is defined as
\begin{eqnarray*}%%%%%%%%%%%%%
[w,u] = \partial_{x_1}^2 w \cdot \partial_{x_2}^2 u+ \partial_{x_2}^2 w \cdot \partial_{x_1}^2 u - 2 \cdot \partial_{x_1 x_2}^2 w \cdot \partial_{x_1 x_2}^2 u, 
\end{eqnarray*}%%----------------------------%%
where $v = v(w)$ is determined by the following problem:
\begin{eqnarray*}%%%%%%%%%%%%%
\Delta v+[w,w] = 0 \text{ in } \Gamma, \quad \partial_\nu v = v =0 \text{ on } \partial\Gamma.
\end{eqnarray*}%%----------------------------%%
This model is used for the plates with large displacements and it is important in post-buckling analysis. More about this model can be found in \cite{karmanplates}.\\
\textbf{The Berger model}: In this case, the nonlinear elastic force takes the form
\begin{eqnarray*}%%%%%%%%%%%%%
\mathcal{F}(w) =-\Big[ \nu \int_{\Gamma} |\nabla w|^2 dx - G  \Big] \Delta w - h,
\end{eqnarray*}%%----------------------------%%
where $\nu >0$ and $G \in \mathbb{R}$ are parameters, $h\in L^2(\Gamma)$. The corresponding potential is 
\begin{eqnarray*}%%%%%%%%%%%%%
\Pi(w) = \frac{\nu}{4} \Big[\int_\Gamma |\nabla w|^2 \Big]^2 - \frac{G}{2}\int_\Gamma |\nabla w|^2  - \int_\Gamma w h.
\end{eqnarray*}%%----------------------------%%
For more details about this model, see \cite[\text{Chapter 7}]{berger}.


\begin{thebibliography}{99}
	
	
	
	\bibitem{adams}
	\newblock R.A. Adams, John J.F. Fournier,
	\newblock {\em Sobolev Spaces}, Second edition,
	\newblock Pure and Applied Mathematics, Vol. 140,
	\newblock Elsevier (Singapore) Ltd., 2003.
	
	
	
	
	
	
	
	
	\bibitem{strongzero}
	\newblock H. Beir\~{a}o da Veiga,
	\newblock On the existence of strong solutions to a coupled fluid-structure evolution problem,
	\newblock {\em J. Math. Fluid Mech.}, 6(2004), 21-52.
	
	\bibitem{bravin}
	\newblock M. Bravin,
	\newblock Energy equality and uniqueness of weak solutions of a “viscous incompressible
	fluid + rigid body” system with Navier slip-with-friction conditions in a 2D bounded
	domain,
	\newblock {\em J. Math. Fluid Mech.}, 21(2019), 23.
	
	
	\bibitem{compressible}
	\newblock D. Breit, S. Schwarzacher,
	\newblock Compressible fluids interacting with a linear-elastic shell, 
	\newblock {\em Arch. Ration. Mech. Anal.}, 228(2017), 495–562.
	
	

	
	
	
	\bibitem{time}
	\newblock A. Chambolle, B. Desjardins, M.J. Esteban, C. Grandmont,
	\newblock Existence of weak solutions for the unsteady interaction of a viscous fluid with an elastic plate,
	\newblock {\em J. Math. Fluid Mech.}, 7(2005), 368-404.
	
	\bibitem{WSUrigid}
	\newblock N. V. Chemetov, \v{S}. Ne\v{c}asov\'{a} , B. Muha,
	\newblock Weak-strong uniqueness for fluid-rigid body interaction problem with slip boundary condition, 
	\newblock {\em J. Math. Phys.}, (2019)60, 011505.
	
	

	
	\bibitem{berger}
	\newblock I. Chueshov, I.Lasiecka,
	\newblock {\em Long-time Behavior of Second Order Evolution Equations with Nonlinear Damping}, 
	Memoirs of AMS, vol 195, AMS, Providence, RI, 2008, No. 912.
	
	\bibitem{karmanplates}
	\newblock I. Chueshov, I. Lasiecka, 
	\newblock {\em Von Karman Evolution Equations}, Springer, New York, 2010.
	
	
	\bibitem{platesproofs}
	\newblock I. Chueshov, S. Kolbasin, 
	\newblock Long-time dynamics in plate models with strong nonlinear damping,
	{\em Comm. Pure App. Anal.}, 11(2012), 659-674.

	
	
	
	
	\bibitem{evans}
	\newblock L. C. Evans,
	\newblock {\em Partial Differential Equations}: Second Edition, 
	Graduate Studies in Mathematics vol. 19, AMS, Providence, Rhode Island, 2010.
	
	
	\bibitem{comWSU1}
	\newblock E. Feireisl, Y. Sun, A. Novotn\'{y},
	\newblock Suitable weak solutions to the Navier–Stokes equations of compressible viscous fluids,
	\newblock {\em Indiana Univ. Math. J.}, 60(2011), 611-631.
	
	\bibitem{comWSU2}
	\newblock E. Feireisl, B. J. Jin, A. Novotn\'{y},
	\newblock Relative entropies, suitable weak solutions, and weak strong uniqueness for the compressible Navier-Stokes system, 
	\newblock {\em J.math. fluid mech.}, 14(2012), 717-730.


	\bibitem{comWSUin}
	\newblock G. P. Galdi, \v{S}. Ne\v{c}asov\'{a}, V. M\'{a}cha,
	\newblock On weak solutions to the problem of a rigid body with a cavity filled with a
	compressible fluid, and their asymptotic behavior,
	\newblock {\em Int. J. Nonlin. Mech.}, 121(2020), 103431.
	
	
	
	
	\bibitem{3dmesh}
	\newblock M. Gali\'{c}, B. Muha, S. \v{C}ani\'{c}
	\newblock Analysis of a 3D nonlinear, moving boundary problem describing fluid-mesh-shell interaction,
	\newblock (2019), https://arxiv.org/abs/1911.09927.
	
		
	
	\bibitem{comWSU0}
	\newblock P. Germain,
	\newblock Weak-strong uniqueness for the isentropic compressible Navier–Stokes system,
	\newblock {\em J. Math. Fluid Mech.}, 13(2011), 137–146.
	
	\bibitem{uniqueness}
	\newblock O. Glass, F. Sueur,
	\newblock Uniqueness results for weak solutions of two-dimensional	fluid-solid systems,
	\newblock {\em Arch. Ration. Mech. Anal.}, 218(2015), 907–944.
	
	\bibitem{grandmont3}
	\newblock C. Grandmont,
	\newblock Existence of weak solutions for the unsteady interaction of a viscous fluid with an elastic plate,
	\newblock {\it SIAM J. Math. Anal.}, 40(2007), 716-737.
	
	\bibitem{grandmont1}
	\newblock C. Grandmont, M. Hillairet,
	\newblock Existence of global strong solutions to a beam-fluid interaction system,
	\newblock {\em Arch. Ration. Mech. Anal.}, 220(2016), 1283-1333.
	
	\bibitem{grandmont2}
	\newblock C. Grandmont, M. Hillairet, J. Lequeurre,
	\newblock Existence of local strong solutions to fluid-beam and fluid-rod interaction systems,
	\newblock {\em Ann. Inst. H. Poincar\'{e} Anal. Non Lin\'{e}aire}, 36(2019), 1105-1149.
	
	\bibitem{continuous}
	\newblock G. Guidoboni, M. Guidorzi, M. Padula,
	\newblock Continuous Dependence on Initial Data in Fluid–Structure Motions,
	\newblock {\em J. Math. Fluid Mech.}, 14(2012), 1–32.
	
	
	
	
	\bibitem{thefirstpaper}
	\newblock D. Hasanyan, N. Hovakimyan, A.J. Sasane, V. Stepanyan, 
	\newblock Analysis of nonlinear thermoelastic plate equations,
	\newblock {\em Proceedings of the 43rd IEEE Conference on Decision and Control}, 2(2004), 1514-1519.
	
	\bibitem{comWSUrigid}
	\newblock O. Kreml, \v{S}. Ne\v{c}asov\'{a}, T. Piasecki,
	\newblock Weak-strong uniqueness for the compressible fluid-rigid body interaction,
	\newblock {\em J. Math. Phys.}, 60(2019), 011505.

	
	
	\bibitem{lasiecka}
	\newblock I. Lasiecka, S. Maad, A. Sasane, 
	\newblock Existence and exponential decay of solutions to a quasilinear thermoelastic plate system,
	{\em Nonlin. Diff. Eq. App.}, 15(2008), 689-715.
	
	
	
	
	\bibitem{ruzicka} 
	\newblock D. Lengeler, M. R\r{u}\v{z}i\v{c}ka, 
	\newblock Weak solutions for an incompressible newtonian fluid interacting with a Koiter type shell,
	\newblock {\em Arch. Ration. Mech. Anal.}, 211(2014), 205–255.
	
	\bibitem{strong}
	\newblock J. Lequeurre,
	\newblock Existence of strong solutions to a fluid-structure system, 
	\newblock {\em SIAM J. Math. Anal.}, 43(2010), 389-410. 
	
	\bibitem{movingdomains}
	\newblock O. Kreml, \v{S}. Ne\v{c}asov\'{a}, T. Piasecki,
	\newblock Local existence of strong solutions and	weak–strong uniqueness for the compressible Navier–Stokes system on moving domains,
	\newblock {\em Proceedings of the Royal Society of Edinburgh Section A: Mathematics}, 150(2020). 2255 - 2300.
	
	
	\bibitem{comstrong}
	\newblock D. Maity, A. Roy, T. Takahashi,
	\newblock Existence of strong solutions for a system of interaction between a compressible viscous fluid and a wave equation,
	\newblock (2020), https://hal-icp.archives-ouvertes.fr/INRIA/hal-02908420v1.
	
	\bibitem{NSFFSI}
	\newblock D. Maity, T. Takahashi,
	\newblock Existence and uniqueness of strong solutions for the system of interaction between a compressible Navier-Stokes-Fourier fluid and a damped plate equation,
	\newblock (2020), https://arxiv.org/abs/2006.00488.
	
	
	
	
	
	\bibitem{sourav}
	\newblock S. Mitra,
	\newblock Local existence of strong solutions for a fluid-structure interaction model, 
	\newblock {\em J. Math. Fluid Mech.}, 22(2018), 60.
	
	
	\bibitem{Boris}
	\newblock B. Muha,
	\newblock A note on the trace Theorem for domains which are locally subgraph of H\"{o}lder continuous function,
	\newblock {\em Networks Hete. Media}, 9(2014), 191-196.
	
	
	
	\bibitem{prodi-serrin}
	\newblock B. Muha, \v{S}. Ne\v{c}asov\'{a}, A. Rado\v{s}evi\'{c}, 
	\newblock A uniqueness result for 3D incompressible fluid-rigid body interaction problem, 
	\newblock {\em J. Math. Fluid Mech.}, 23(2020), 1.
	
	
	
	\bibitem{muhasch}
	\newblock B. Muha, S. Schwarzacher,
	\newblock Existence and regularity for weak solutions for a fluid interacting with a non-linear shell in 3D,
	\newblock (2019), https://arxiv.org/abs/1906.01962.
	
	
	
	\bibitem{BorSunNavierSlip}
	\newblock B. Muha, S. {\v C}ani{\'c},
	\newblock Existence of a weak solution to a fluid-elastic structure interaction problem with the Navier slip boundary condition,
	\newblock {\em J. Diff. Equ.}, 260(2016), 8550-8589.
	
	
	
	\bibitem{BorSun}
	\newblock B. Muha, S. {\v C}ani{\'c},
	\newblock Existence of a weak solution to a nonlinear
	fluid-structure interaction problem modeling the flow of an
	incompressible, viscous fluid in a cylinder with deformable
	walls,
	\newblock {\em Arch. Ration. Mech. Anal.}, 207(2013),  919-968.
	
	
	
	
	
	\bibitem{WSinc}
	\newblock S. Schwarzacher, M. Sroczinski,
	\newblock Weak-strong uniqueness for an elastic plate interacting with the Navier Stokes equation,
	\newblock (2020), https://arxiv.org/abs/2003.04049.
	
	\bibitem{triebel}
	\newblock H. Triebel,
	\newblock Theory of function spaces II, volume 84 of Monographs in Mathematics,
	\newblock Birkh\:{a}user Verlag, Basel, 1992.
	
	\bibitem{trwa}
	\newblock S. Trifunovi\'c, Y.-G. Wang, 
	\newblock Existence of a weak solution to the fluid-structure interaction problem in 3D,
	\newblock {\em J. Diff. Equ.} 268 (2020), 1495-1531.
	

	
	\bibitem{trwa3}
	\newblock S. Trifunovi\'c, Y.-G. Wang, 
	\newblock On the interaction problem between a compressible viscous fluid and a nonlinear thermoelastic plate,
	\newblock (2020), https://arxiv.org/abs/2010.01639.
	
	
	
	
	
	
	
	
	
\end{thebibliography}
\end{document}